\documentclass[reqno, a4paper, 11pt]{amsart}
\usepackage{amsmath, amsfonts, amsthm, amssymb}
\numberwithin{equation}{section}

\newcommand{\Q}{\mathbb{Q}}

\newcommand{\R}{\mathbb{R}}
\newcommand{\C}{\mathbb{C}}
\newcommand{\N}{\mathbb{N}}
\newcommand{\Z}{\mathbb{Z}}

\newtheorem{thm}{Theorem}
\newtheorem{lem}{Lemma}
\newtheorem{pro}{Proposition}

\renewcommand{\mod}[1]{\hspace{-2.9mm}\pmod{#1}}
\newcommand{\x}{{\bf x}}
\newcommand{\e}{e}
\newcommand{\rom}{\mathrm}
\newcommand{\bfP}{\mathbb{P}}

\newcommand{\ov}{\overline}
\newcommand{\ma}{\mathbf}
\newcommand{\ben}{\begin{enumerate}}
\newcommand{\een}{\end{enumerate}}
\newcommand{\eit}{\begin{itemize}}

\newcommand{\ve}{\varepsilon}
\newcommand{\de}{\delta}
\newcommand{\mcal}{\mathcal}
\newcommand{\lab}{\label}
\newcommand{\al}{\alpha}

\newcommand{\D}{\Delta}
\newcommand{\be}{\beta}
\newcommand{\la}{\lambda}

\newcommand{\colt}[2]{\genfrac{}{}{0pt}{1}{#1}{#2}}
\newcommand{\eqm}[2]{\equiv #1 \pmod{#2}}

\newcommand{\vr}{\varrho}
\newcommand{\vt}{\vartheta}
\renewcommand{\d}{\mathrm{d}}
\renewcommand{\le}{\leqslant}
\renewcommand{\leq}{\leqslant}
\renewcommand{\geq}{\geqslant}
\newcommand{\bxi}{\boldsymbol{\xi}}
\newcommand{\bt}{\boldsymbol{\tau}}
\newcommand{\bla}{\boldsymbol{\lambda}}
\newcommand{\m}{\mathfrak{m}}
\newcommand{\E}{\mathcal{E}}
\newcommand{\FF}{\mathcal{F}}

\newcommand{\tstack}[3]{\colt{#1}{\colt{#2}{#3}}}

\newcommand\HH{\mathcal{H}}

\newcommand\f{{f}}

\newcommand\NN{\mathcal{N}}
\newcommand\cR{\mathcal{R}}
\newcommand\TT{\mathcal{T}}

\newcommand\PP{\mathbb{P}}
\newcommand\xx{\mathbf{x}}
\newcommand\dd{\,\mathrm{d}}
\newcommand\Pthree{{\PP^3}}
\newcommand{\base}[7]{\xi^{({#1},{#2},{#3},{#4},{#5},{#6},{#7})}}
\newcommand{\congr}[3]{{#1} \equiv {#2}\ (\mathrm{mod}\ {#3})}

\newcommand\NU{{N_{U,H}}}
\newcommand{\Esix}{{\mathbf E}_6}
\newcommand{\Dfour}{{\mathbf D}_4}

\newcommand{\tS}{{\widetilde S}}
\newcommand{\anti}{-K_\tS}

\DeclareMathOperator{\Pic}{Pic}
\DeclareMathOperator{\Spec}{Spec}

\DeclareMathOperator{\hcf}{gcd}

\DeclareMathOperator{\rk}{rk}

\theoremstyle{definition}

\newtheorem*{ack}{Acknowledgements}

\begin{document}

\title[On Manin's conjecture for a cubic surface]{On Manin's
  conjecture for a certain singular cubic surface}

\author{R. de la Bret\`eche}
\author{T.D. Browning}
\author{U. Derenthal}

\address{
Institut de Math\'ematiques de Jussieu,
Universit\'e Paris 7 Denis Diderot,
Case Postale 7012,
2, Place Jussieu, 
F-75251 Paris cedex 05}
\email{breteche@math.jussieu.fr}

\address{School of Mathematics,
University of Bristol, Bristol BS8 1TW}
\email{t.d.browning@bristol.ac.uk}

\address{Mathematisches Institut,
Bunsenstr. 3--5, D-37073 G\"ottingen}
\email{derentha@math.uni-goettingen.de}


\begin{abstract}
This paper contains a proof of the Manin conjecture for
the singular cubic surface $S \subset \Pthree$ that is defined by the equation
$x_1x_2^2+x_2x_0^2+x_3^3=0$. In fact if $U \subset S$
is the Zariski open subset obtained by deleting the unique line from
$S$,  and $H$ is the usual exponential height on
$\bfP^3(\Q)$, then the height zeta function
$
\sum_{x \in U(\Q)}{H(x)^{-s}}
$
is analytically continued to the half-plane $\Re e (s)>9/10$.
\newline
\newline
\noindent
%{\em Sur la conjecture de Manin pour une certaine
%surface cubique singuli\`ere.}\\
\textsc{R\'esum\'e.}
Ce papier contient une preuve de la conjecture de Manin pour la surface
cubique singuli\`ere $S \subset \Pthree$ d\'efinie par
$x_1x_2^2+x_2x_0^2+x_3^3=0$.
En effet, si $U \subset S$ est l'ouvert obtenu en enlevant l'unique
droite contenue dans $S$  et $H$ est la fonction des hauteurs usuelle
de $\bfP^3(\Q)$,
alors la fonction z\^eta des hauteurs
$
\sum_{x \in U(\Q)}{H(x)^{-s}}
$
peut \^etre prolong\'ee de mani\`ere analytique au demi-plan $\Re e
(s)>9/10$.
\end{abstract}

\maketitle

\section{Introduction}\lab{intro}

Let $S \subset \bfP^3$ be a cubic surface that is
defined over $\Q$ and has isolated singularities. As soon as $S$
contains a single $\Q$-rational point the set of
rational points $S(\Q)=S\cap\bfP^3(\Q)$ is dense in the Zariski
topology, and it is natural to seek a finer interpretation of this density.
Given a point $x=[x_0,\ldots,x_3] \in \bfP^3(\Q),$ with
$x_0,\ldots,x_3 \in \Z$ such that $\hcf(x_0,\ldots,x_3)=1$, we let
$$
H(x)=\max\{|x_0|,|x_1|,|x_2|,|x_3|\}.
$$
Then $H: \bfP^3(\Q) \rightarrow \R_{\geq 0}$ is
the exponential height attached to
the anticanonical embedding of $S$, metrized by the choice of norm
$|\ma{z}|:=\max_{0\leq i \leq 3}|z_i|$ on
$\R^{4}$.  We may define the quantity
$$
N_{U,H}(B)=\#\{x \in U(\Q): H(x) \leq B\},
$$
for any $B\geq 1$ and any Zariski open subset $U\subseteq S$.
If $S$ contains lines defined over $\Q$ then $N_{S,H}(B)$ will be
dominated by the rational points of height at most $B$ that lies on
such lines.  For this reason one is most interested in studying the
counting function $N_{U,H}(B)$ for the open subset $U\subset S$
obtained  by deleting all of the lines from $S$.

In this setting Manin \cite{f-m-t} has formulated a
far-reaching conjecture for the asymptotic
behaviour of $N_{U,H}(B)$, as $B \rightarrow
\infty$.
This states that there is a non-negative constant
$c_{S,H}$ and a positive integer
$\rho$ such that
\begin{equation}\lab{manin}
\NU(B) = c_{S,H} B (\log B)^{\rho-1}\big(1+o(1)\big),
\end{equation}
as $B\rightarrow \infty$.  Here $\rho$ is conjectured to be the
rank of the Picard group of the minimal desingularisation of $S$,
and the constant $c_{S,H}$ has also been given a conjectural interpretation
at the hands of Peyre \cite{MR1340296}, Batyrev and
Tschinkel~\cite{b-t}, and Salberger~\cite{MR1679841}.

Although Manin's conjecture can actually be applied to a rather
general class of algebraic variety, in which context it has met with a
reasonable degree of success, the
situation for cubic surfaces is rather less satisfactory.
For non-singular cubic surfaces the best result that we have is the
upper bound  $\NU(B)=O_{\varepsilon,S}(B^{4/3+\varepsilon})$,
which is due to Heath-Brown \cite{MR98h:11083} and
applies when the surface contains three coplanar lines defined over $\Q$.
For singular cubic surfaces better estimates are available.
A modern classification of such surfaces can be found in the work of Bruce and Wall \cite{MR80f:14021},
which shows in particular that there are only finitely
many isomorphism classes to consider over $\ov{\Q}$, 
these being essentially classified by their singularity type.
The Manin conjecture for singular cubic surfaces has only been settled
in particularly simple cases, such as the singular toric variety
$$
x_0^3=x_1x_2x_3
$$
of singularity type $3\textbf{A}_2$.
Several authors have studied this surface, and the sharpest
estimate available is that due to the first author \cite{MR2000b:11074}.
Further work worth mentioning is that due to Chambert-Loir and Tschinkel
\cite{ct}, who have established Manin's conjecture for any cubic surface
arising as an equivariant compactifications of $\mathbb{G}_a^2$. There
is also the work of Heath-Brown \cite{0981.32025} and the second author \cite{math.NT/0404245}.
These latter results provide upper and lower bounds
of the expected order of magnitude for the counting function
associated to two singular cubic surfaces: the Cayley cubic surface
$$
\frac{1}{x_0}+\frac{1}{x_1}+\frac{1}{x_2}+\frac{1}{x_3}=0
$$
of singularity type $4\mathbf{A}_1$, and the surface
$$
x_1x_2x_3=x_0(x_1+x_2+x_3)^2
$$
containing a $\Dfour$-singularity, respectively.

We are now ready to reveal the contribution that we have been able to
make to this topic.  Thus the primary goal of this paper is to verify
the Manin conjecture for the cubic surface
\begin{equation}\label{eq:surface}
x_1x_2^2+x_2x_0^2+x_3^3 = 0,
\end{equation}
which we henceforth denote by $S$.  This surface contains a unique
singularity of type $\Esix$, and a unique line $\ell$ which is given by
$x_2=x_3 = 0$.  It has been shown by Hassett and Tschinkel \cite[Remark 4.3]{MR2029868} that $S$
is not an equivariant compactification of $\mathbb{G}_a^2$, so that it
is not covered by \cite{ct}. Let $U\subset S$  be the open subset formed by deleting
$\ell$ from $S$. Then we have the following result.

\begin{thm}\lab{main'}
Let $\varepsilon>0$.  Then there exists a polynomial $P$ of
degree $6$ such that
$$
\NU(B)=B P(\log B) +O_\varepsilon(B^{10/11+\varepsilon}),
$$
for any $B \geq 1$.
Moreover the leading coefficient of $P$ is equal to
$$
\frac{\omega_\infty}{6220800}
\prod_p\Big(1-\frac{1}{p}\Big)^7\Big(1+\frac 7 p + \frac 1 {p^2}\Big),
$$
where
\begin{equation}\lab{om-inf}
\omega_\infty = 12 \int\int\int_{\{(t,u,v)\in
   \R^3 : ~|t^2+u^3|\le 1, ~0\leq tv^3 \leq 1, ~0
\le  v \le 1, ~|uv^4| \le 1\}} \dd t \dd u \dd v.
\end{equation}
\end{thm}

We shall verify in \S \ref{sec:manins_conjecture} that Theorem
\ref{main'} is in agreement with the Manin conjecture.
In the classification of singular cubic surfaces over $\ov{\Q}$, $S$
the only cubic surface with an $\Esix$-singularity, up to projectivity
\cite{MR80f:14021}. 
In fact this is the most extreme type of  singularity that a cubic surface
can possess.
Given that non-singular cubic surfaces seem so difficult
to tackle, our success with \eqref{eq:surface} perhaps reflects the
fact that we are as far away from the non-singular setting as possible.

It is now well-recognised that universal torsors play a central
r\^ole in proofs of  the Manin conjecture for Fano varieties.
There is no exception to this philosophy in the present work. Thus
in the proof of Theorem~\ref{main'} crucial use is made of the
universal torsor above the minimal desingularisation $\tS$ of $S$,
which turns out to have the natural affine embedding
\begin{equation}\lab{ut-intro}
\tau_\ell\xi_\ell^3\xi_4^2\xi_5+\tau_2^2\xi_2+\tau_1^3\xi_1^2\xi_3 = 0.
\end{equation}
This has been calculated by Hassett and Tschinkel using the Cox ring \cite{MR2029868}.
However in \S \ref{proof2}  we shall establish a completely explicit bijection between
$U(\Q)$ and a suitable set of integral points satisfying this equation,
via an elementary analysis of the equation \eqref{eq:surface} defining $S$.
It will become apparent that the passage to the universal torsor is really
only the first step on the road to proving Theorem~\ref{main'}, and
that a considerable amount of input is still required.

Once the passage to the universal torsor is accomplished, the proof of Theorem
\ref{main'} broadly follows the strategy of the first two authors
\cite{math.NT/0412086, math.NT/0502510},
where key use is made of the fact that the torsor equation in each
case has precisely three terms.
In counting integral solutions to \eqref{ut-intro}, subject to certain
constraints, we shall begin by fixing most of the variables and
summing only over the variables
$\tau_1,\tau_2,\tau_\ell$.
The key idea is then to view the  equation as a congruence
$$
\tau_2^2\xi_2\equiv -\tau_1^3\xi_1^2\xi_3 ~~\mod{\xi_\ell^3\xi_4^2\xi_5},
$$
in order to take care of the summation over $\tau_\ell$.
One proceeds to employ standard facts about the number of integer solutions to
polynomial congruences that are restricted to lie in certain regions.
This produces a main term and an error term, and the rest of the
proof involves summing each of these terms over all of the
remaining variables.  While the treatment of the main term is
relatively routine, the treatment of the error term presents a much
more serious obstacle.  There are two main ingredients in this part of
the work, both of which are rooted in the theory of exponential sums.
The first involves showing that sequences of the form
$(ax^3+bx^2)/q$ are equidistributed modulo $1$ as $x$ ranges over
the ring $\Z/q\Z$, for fixed integers $a,b,q$ such that
$\hcf(a,b,q)=1$, and the second constitutes a
delicate analysis of certain exponential sums involving
real-valued functions that arise in our work.
Whereas the first ingredient is independent of the choice of
norm used to metrize the height function $H$, and so may be
thought of as purely ``arithmetic'', the second ingredient is
intimately connected to the norm selected and may be thought of as
being ``analytic'' in nature.

Given the shape of the estimate in Theorem \ref{main'} it is no
surprise that we are able to say something about the corresponding
height zeta function. As above let $U\subset S$ be the open subset of the surface
\eqref{eq:surface} that is formed by deleting the
unique line from it. Then we may define
$$
Z_{U,H}(s):=\sum_{x \in U(\Q)}\frac{1}{H(x)^s},
$$
for $\Re e(s)>1$, and Theorem \ref{main'} can be used to show that
$Z_{U,H}(s)$ has a meromorphic continuation to the half-plane $\Re
e(s)>10/11$.  In fact by returning to the proof of Theorem~\ref{main'}
we are able to say something about the analytic structure of
$Z_{U,H}(s)$ to the left of the line $\Re e(s)=10/11$.
For $\Re e (s)>0$ we define the functions
\begin{align}
E_1(s+1)
&:=\zeta(2s+1)\zeta(3s+1)^2\zeta(4s+1)^2\zeta(5s+1)\zeta(6s+1),
\lab{E1}
\\
E_2(s+1)
&:=\frac{\zeta(13s+3)^5\zeta(14s+3)^2}{\zeta(7s+2)^4\zeta(8s+2)^4\zeta(9s+2)^2\zeta(10s+2)\zeta(19s+4)^2}. 
\lab{E2}
\end{align}
It is easily seen that $E_1(s)$ has a meromorphic
continuation to the entire complex plane with a single pole at $s=1$,
and similarly,  $E_2(s)$ is holomorphic and bounded
on the half-plane $\Re e (s)>9/10$.
We are now ready to record precisely what we have been able to say
about the height zeta function.

\begin{thm}\lab{thm:main}
Let $\varepsilon>0$.  Then there exists a constant $\beta
\in \R$, and functions $G_1(s), G_2(s)$
that are holomorphic on the half-plane $\Re
e(s)\geq 43/48+\varepsilon$, such that for $\Re e(s)>1$ we have
$$
Z_{U,H}(s) =E_1(s)E_2(s)G_1(s) +\frac{{12/\pi^2}+2\be}{s-1}+G_2(s).
$$
In particular $(s-1)^7 Z_{U,H}(s)$ has a holomorphic
continuation to the half-plane $\Re e(s)>9/10.$
\end{thm}

Explicit expressions for $\beta, G_1$ and $G_2$ can be found in
\eqref{defbeta}, \eqref{G1} and \eqref{G2}, respectively.
It can be seen there that $G_1(s)$ is actually holomorphic and bounded
on the half-plane $\Re e(s)\geq 5/6+\varepsilon$, and that
$$
G_2(s)\ll
1+|\Im m (s)|
$$
for $\Re e(s)\geq 43/48+\varepsilon$.

With more work it is likely that the constant $43/48$ can be
reduced slightly, although all we need to deduce the final sentence in
Theorem \ref{thm:main} is the fact that $43/48<9/10$.
However, under the assumption of the Riemann hypothesis it is clear that $E_2(s)$ is holomorphic for $\Re e
(s)>8/9$,  whence $Z_{U,H}(s)$ has a meromorphic
continuation to the half-plane $\Re e(s)>43/48.$

Theorem \ref{thm:main} bears a striking resemblance to the results
obtained by the first two authors \cite{math.NT/0412086,
math.NT/0502510}, in their work on the Manin conjecture for singular
del Pezzo surfaces of degree 4, which
also contain explicit expressions for the
corresponding height zeta functions.
Thus in addition to the ``main term'' $E_1(s)E_2(s)G_1(s)$, all of
these results have a term $\frac{12}{\pi^2}(s-1)^{-1}$
that corresponds here to the residual conic obtained by intersecting $S$
with the plane $x_3=0$, and a further ``$\beta$-term''.
In Theorem \ref{thm:main} the constant $\beta$ has much in common
with the corresponding result in \cite{math.NT/0412086}, arising as it
does through  the application of results about
the equidistribution of squares in a fixed residue class.
However the argument needed here is distinctly subtler
than anything previously encountered.

The genesis of this paper lies in an earlier paper due to the third
author \cite{ulrich}, who succeeded in proving a version of Theorem \ref{main'} with
an error term of $O(B(\log B)^2)$.  The main contribution of the
first and second author has therefore been to push the analysis
further, to the extent that we now have results of the
precision detailed above.  During the final preparation of this paper, the authors have been made
aware of the doctoral thesis of M. Joyce at Brown University, who has
independently established the Manin conjecture
for the $\Esix$ cubic surface~$S$. His main result is weaker than that obtained in our
paper, since he only establishes an asymptotic
formula with an error term of  $O(B (\log B)^5)$.

We end this introduction by giving an overview of the contents of this paper.
As indicated above, we shall begin in \S \ref{sec:manins_conjecture}
by showing that Theorem~\ref{main'} is in
complete agreement with the Manin conjecture.
Next in \S \ref{e-sum} and \S \ref{equi} we shall collect together
most of the material concerning exponential sums and equidistribution
that will be crucial for our treatment of the error terms discussed above.
In \S \ref{real}
we shall introduce and analyse a number of real-valued functions that
will arise in our work, before turning in \S \ref{proof1} to a
preliminary estimate for the counting function $\NU(B)$.  The passage to
the universal torsor will take place in \S \ref{proof2}, and the
conclusion of the proof of Theorem \ref{main'} will form the contents
of \S \ref{proof3} and \S \ref{proof4'}.
Finally we shall deduce the statement of Theorem \ref{thm:main} in \S
\ref{proof4}.

\begin{ack}
Part of this work was undertaken while the second author was attending
the \emph{Diophantine Geometry} intensive research period at the
Centro di Ricerca Matematica Ennio De Giorgi in Pisa, and the third
author was visiting Brendan Hassett at Rice University.
The paper was finalised while the first author was
at the \'Ecole Normale Sup\'erieure,
and the second author was at Oxford University supported by
EPSRC grant number GR/R93155/01.
The hospitality and financial support of all these institutions is
gratefully acknowledged.
Finally, the authors would like to thank the anonymous
referee for his extremely attentive reading of the manuscript and
numerous helpful suggestions. In particular, the referee's comments
have helped to simplify the proofs of
Lemmas \ref{exp1} and \ref{exp3}, and led to an overall improvement in Lemma \ref{except}.
\end{ack}

\section{Conformity with the Manin conjecture}\label{sec:manins_conjecture}

In this section we shall review some of the geometry of the surface $S\subset
\bfP^3$, with a view to calculating the invariants
appearing in Manin's conjecture and its refinement by Peyre.
Let $\tilde S$ denote the minimal desingularisation of $S$,  and let
$\pi: \tS\rightarrow S$ denote the corresponding
blow-up map.    We let $F_1,\ldots,F_6$ denote the exceptional divisors
of $\pi$.  Then the divisors $F_1,\ldots,F_6$ are all defined over $\Q$, and
together with the line $\ell$,
they generate the Picard group $\Pic(\tS)$ of $\tS$.
In particular we have $\rho=7$ in \eqref{manin}.

Turning to the conjectured value of the constant $c_{S,H}$ in
\eqref{manin}, we follow the notation and methodology of Peyre \cite{MR1340296}.
With this in mind we proceed by establishing the following result.

\begin{lem}\label{lem:manin_conjecture}
We have  $c_{S,H} =\alpha(\tS)\beta(\tS)\omega_H(\tS)$, with
$$
\alpha(\tS)= \frac 1{6220800}, \quad
\beta(\tS) = 1,
\quad
\omega_H(\tS) = \omega_\infty
\prod_p\Big(1-\frac{1}{p}\Big)^7\Big(1+\frac 7 p
+ \frac 1 {p^2}\Big),
$$
where $\omega_\infty$ is given by \eqref{om-inf}.
\end{lem}

\begin{proof}
We have already observed that $\{F_1,F_2,F_3,\ell,F_4,F_5,F_6\}$ is a
basis for
$\Pic(\tS)$.  It follows from \cite{MR2029868} that the effective cone
$\Lambda_\rom{eff}(\tS)$ is equal to $\Pic(\tS) \otimes_\Z \R$,
and that the  dual cone of nef divisors is simplical, in the sense
that it is generated by $\rho=7$ elements.
Moreover
the anticanonical divisor $\anti$ of $\tS$ is given by
\[
\anti = 2F_1 + 3F_2 + 4F_3+3\ell+4F_4+5F_5+6F_6.
\]
We may therefore write $\anti=\bla$ in the basis
$\{F_1,F_2,F_3,\ell,F_4,F_5,F_6\}$, with
\begin{equation}\lab{bla}
\bla=(\la_1,\la_2,\la_3,\la_\ell,\la_4,\la_5,\la_6):=(2,3,4,3,4,5,6).
\end{equation}
Thus the definition of $\alpha(\tS)$ reveals that
\begin{align*}
\alpha(\tS)
&= \rom{meas}\{ \ma{t}\in \R_{\geq 0}^7: ~\bla. \ma{t}=1 \}
= \frac{1}{6!\la_1\la_2\la_3\la_\ell\la_4\la_5\la_6}=\frac{1}{6220800},
\end{align*}
where we have written $\ma{t}=(t_1,t_2,t_3,t_\ell,t_4,t_5,t_6)$.
Next we note that
\[
\beta(\tS) :=
\#H^1(\rom{Gal}(\overline{\Q}/\Q),\Pic(\tS)\otimes_\Q
\overline{\Q})= 1,
\]
since $S$ is split over $\Q$.
Finally we must consider the factor $\omega_H(\tS)$, which corresponds
to a product of local densities.  According to the definition of
$\omega_H(\tS)$ we have
\begin{align*}
\omega_H(\tS)
&:= \lim_{s \to 1}((s-1)^{\rk
      \Pic(\tS)} L(s,\Pic(\tS))) \omega_\infty \prod_p \frac{\omega_p}{L_p(1,\Pic(\tS))}\\
&= \omega_\infty \prod_p \Big(1 - \frac{1}{p}\Big)^7 \omega_p,
\end{align*}
since $L(s, \Pic(\tS)) = \zeta(s)^7$,
in our case.  The calculation of $\omega_p$ is straightforward, and
ultimately leads to the conclusion that
$$
\omega_p= 1+\frac{7}{p}+\frac{1}{p^2}.
$$
To compute $\omega_\infty$ we parametrise the points by writing $x_1$ as a function of
$x_0,x_2,x_3$ in $f(\x)=x_1x_2^2+x_2x_0^2+x_3^3$.
Since $\xx = -\xx$ in $\Pthree$, we may assume $x_2
\ge 0$. On observing that
$\frac{\partial\f}{\partial x_1} = x_2^2$, the
Leray form
$\omega_L(\tS)$ is given by $x_2^{-2} \dd x_0 \dd x_2 \dd x_3$,
and so
\[
\omega_\infty=
2\int\int\int_{\{|x_2^{-2}(x_2x_0^2+x_3^3)|\le 1, ~0 \le x_0, x_2
\le 1, ~|x_3|\le 1\}} x_2^{-2} \dd x_0 \dd x_2 \dd x_3.
\]
But then the change of variables
$x_0 = tx_2^{1/2}$, $x_3 = ux_2^{2/3}$ and $x_2 = v^6$, easily yields
the value of $\omega_\infty$ given in \eqref{om-inf}.
This complete the proof of the lemma.
\end{proof}

On combining Lemma \ref{lem:manin_conjecture} with our earlier
observation that $\rho=7$ in \eqref{manin}, we therefore conclude that
Theorem \ref{main'} is in accordance with the Manin conjecture.

\section{Exponential sums}\lab{e-sum}

During the course of the subsequent section we shall need good upper
bounds for the modulus of several exponential
sums.  We have collected together the
results that we shall need in the present section, throughout which we
employ the usual notation $\e(x)={\e}^{2\pi i x}$ and $\e_q(x)={\e}(x/q)$, for
any $q \in \N$ and $x \in \R$, and always
take $\N$ to denote the set of positive integers.
Furthermore, we shall write $\lfloor x \rfloor$
(resp. $\lceil x \rceil$) for the integer part (resp. the ceiling) of
any $x\in \R$.

Let $a,b \in \Z$ and let $q \in \N$.
The primary goal of this section is then to estimate the cubic exponential sum
\begin{equation}\lab{exp-sum}
S_q(a,b):=\sum_{\colt{x=1}{\hcf(x,q)=1}}^q \e_q(ax^3+bx^2),
\end{equation}
under the assumption that $\hcf(a,b,q)=1$.
Our approach will involve relating $S_q(a,b)$ to the complete exponential sum
\begin{equation}\lab{exp-sum:complete}
T_q(a,b):=\sum_{x=1}^q \e_q(ax^3+bx^2).
\end{equation}
We begin by recording the multiplicativity properties
\begin{equation}\lab{t-m}
\begin{split}
S_{uv}(a,b)&=S_u(v^2a,vb)S_v(u^2a,ub),\\
T_{uv}(a,b)&=T_u(v^2a,vb)T_v(u^2a,ub),
\end{split}
\end{equation}
that are valid for any coprime $u,v \in \N$ such that
$\hcf(a,b,uv)=1$.  These equalities follows from the Chinese remainder
theorem (see \cite[Lemma 2.10]{vaughan}, for example).
We are now ready to estimate \eqref{exp-sum} in the case $b=0$.

\begin{lem}\lab{exp1}
Let $\varepsilon>0$ and suppose that $\hcf(a,q)=1$.  Then we have
$$
S_q(a,0) \ll_\varepsilon q^{2/3+\varepsilon}.
$$
\end{lem}

\begin{proof}
In view of \eqref{t-m} and the estimate $A^{\omega(q)}=O_{A,\ve}(q^\ve)$, it
will suffice to show that $S_{p^\ell}(a,0)\ll p^{2\ell/3}$, for any
prime $p$ such that $p\nmid a$, and any $\ell\in \N$.
But when $\ell \geq 3$ it follows that
$$
S_{p^\ell}(a,0)=T_{p^\ell}(a,0)-p^2T_{p^{\ell-3}}(a,0),
$$
whence \cite[Eq. (7.9)]{vaughan} yields
$$
S_{p^\ell}(a,0)\ll p^{2\ell/3}+p^2p^{2(\ell-3)/3}\ll p^{2\ell/3},
$$
when $\ell \geq 3$. The same sort of calculation suffices to handle
the cases $\ell=1$ and $\ell=2$, which therefore completes the proof
of the lemma.
\end{proof}

We now turn to the task of estimating
\eqref{exp-sum} for non-zero values of $b$, for which
we shall need a corresponding
estimate for \eqref{exp-sum:complete} in the case that $b$ is
non-zero.  This is provided for us by the
following result.

%% see e6-L3.tex for full proof of this result

\begin{lem}\lab{exp2}
Let $p$ be a prime such that $\hcf(a,b,p)=1$ and let $\ell \in \N$. Then we have
$$
|T_{p^\ell}(a,b)| \leq 2 p^{\ell/2}\hcf(b,p^\ell).
$$
\end{lem}

\begin{proof}
The case in which $\ell=1$ is handled by the well-known estimate
of Weil \cite{weil}, which gives $|T_p(a,b)|\leq 2p^{1/2}$.
The case in which $\ell \geq 2$ follows from the work of Loxton and
Vaughan \cite[Theorem 1]{lox}.  This completes the proof of Lemma \ref{exp2}.
\end{proof}

We are now ready to record an estimate for \eqref{exp-sum} that is
valid for any choice of $a,b\in \Z$ and $q \in \N$ such that
$\hcf(a,b,q)=1$.

\begin{lem}\lab{exp3}
Let $\varepsilon>0$ and suppose that $\hcf(a,b,q)=1$.  Then we have
$$
S_q(a,b) \ll_\varepsilon q^{1/2+\varepsilon}\hcf(b,q).
$$
\end{lem}

\begin{proof}
As in the proof of Lemma \ref{exp1}, the properties in
\eqref{t-m} render it sufficient to establish the
bound $S_{p^\ell}(a,b)\ll p^{\ell/2}\hcf(b,p^\ell)$, for any
prime $p$ such that $p\nmid \hcf(a,b)$, and any $\ell\in \N$.
When $\ell \geq 2$ it follows that
$$
S_{p^\ell}(a,b)=T_{p^\ell}(a,b)-pT_{p^{\ell-2}}(ap,b),
$$
whence Lemma \ref{exp2} yields
$S_{p^\ell}(a,b)\ll p^{\ell/2}$, if $\ell \geq 2$ and $p\nmid b$. If $p\mid
b$, then we may write $b=pb'$. In this case Lemma \ref{exp2} yields
\begin{align*}
S_{p^\ell}(a,b)&=T_{p^\ell}(a,b)-p^2T_{p^{\ell-3}}(a,b') \\
&\ll
p^{\ell/2}\hcf(b,p^\ell)+
p^{2+(\ell-3)/2}\hcf(b',p^{\ell-3})\\
&\ll p^{\ell/2}\hcf(b,p^\ell),
\end{align*}
if $\ell \geq 3$.
Together these two estimates handle the case in
which $\ell \geq 3$. Finally, the same sort of calculation suffices to handle
the cases $\ell=1$ and $\ell=2$, which therefore completes the proof
of Lemma \ref{exp3}.
\end{proof}

Now let $I=[t_1,t_2]\subset \R$ be any closed interval, and let
$f$ be a real-valued function on it.  Then for given
$a,b ,q\in \Z$ such that $q>0$, the remainder of this section is
concerned with the size of the exponential sum
\begin{equation}\lab{exp-sum''}
A_{I}(q;a,b,f) :=\sum_{t_1< n\leq t_2} \e_q(a n+b f(n)).
\end{equation}
In particular we shall want to obtain a
saving over the trivial upper bound
\begin{equation}\lab{triv1}
|A_{I}(q;a,b,f)| \leq t_2-t_1+1,
\end{equation}
by restricting our attention to suitable families of real-valued functions.
For an interval $I=[t_1,t_2]\subset \R$ and a real number $\la_0\geq
1$, we shall say that a real-valued function $f$ belongs to the set
$C^1(I;\la_0)=C^1(t_1,t_2;\la_0)$ if $f$ is differentiable on $I$,
with
\begin{equation}\lab{def-sig}
|f(t_2)-f(t_1)|+1\leq \la_0, %t%
\end{equation}
%t%for some absolute positive constant $A$,
and if $f'$ is monotonic and of constant sign on $(t_1,t_2)$.
We then  have the following result.

\begin{lem}\lab{triv2}
Let $I \subset \R$ be any  closed interval and let $\la_0 \geq 1$.
Suppose that $a,b,q\in \Z$ such that $0<|a|\leq q/2$, and
let $f \in C^1(I;\la_0)$.
Then we have
$$
A_{I}(q;a,b,f) \ll
\frac{1}{|a|}\Big(q+|b|\la_0 \Big).
$$
\end{lem}

\begin{proof}
Suppose that $I=[t_1,t_2]$, for $t_1<t_2$.
To establish the lemma, we write $A_t(q;a)$ for the linear exponential sum
$A_{[t_1,t]}(q;a,0,0)$ for $t\in (t_1,t_2]$.  Then
\begin{equation}\lab{p1}
A_t(q;a) = \frac{\e_q(a\lceil t_1 \rceil
)-\e_q(a(\lfloor t \rfloor+1))}{1-e(a/q)},
\end{equation}
whence
\begin{equation}\lab{p2}
A_t(q;a) \ll \frac{1}{|1-{\e}(a/q)|} =\frac{1}{|\sin(\pi a/q)|}\ll
\frac{q}{|a|},
\end{equation}
since $|a|\leq q/2$.
Set $F(t)=\e_q(b f(t))$ for $t_1<t\leq t_2$, and
$F(t)=0$ otherwise.   Then in view of \eqref{p1} and \eqref{p2},
a simple application of partial summation yields
\begin{equation}\lab{pisa2}
\begin{split}
A_{I}(q;a,b,f)&=A_{t_2}(q;a)F(t_2)-\int_{t_1}^{t_2}A_{t}(q;a)F'(t)\d t\\
&=-\int_{t_1}^{t_2}A_{t}(q;a)F'(t)\d t+ O(|a|^{-1}q)\\
&= \int_{t_1}^{t_2} \frac{\e_q(a(\lfloor t\rfloor+1))}{1-{\e}(a/q)}F'(t)\d t +
O(|a|^{-1}q).
\end{split}
\end{equation}
But then the lemma easily follows from the observation that
$$
\int_{t_1}^{t_2} |F'(t)|\d t \leq \frac{2\pi |b|}{q}
\int_{t_1}^{t_2} |f'(t)F(t)|\d t,
$$
this latter integral being $O(\la_0)$.
\end{proof}

We can do somewhat better by further restricting the class of
functions $f$ under consideration. Let $I\subset \R$ be a closed interval, and let
$j,\la_0,\la_1,\la_2 \in \R$ such that
\begin{equation}\lab{la_}
j,\la_0,\la_1\geq 1, \quad \la_2 >0.
\end{equation}
We say that a real-valued function $f$ belongs to the set
$C^2(I;\la_0,\la_1,\la_2,j)$ if $f$
is twice differentiable on $I$, with $f \in C^1(I;\la_0)$ and
$$%t%
|f'(t)|\leq  \la_1,
\quad
\la_2\leq |f''(t)| \leq j\la_2,
$$
throughout $I$.  On defining the notation
\begin{equation}\lab{meas}
\m(I):=\rom{meas}(I)+1,
\end{equation}
we then have the following result.

\begin{lem}\lab{triv2'}
Let $I \subset \R$ be any closed interval and let
$j,\la_0,\la_1,\la_2 \in \R$ such that \eqref{la_}
holds.  Suppose that $a,b,q\in \Z$ such that $0<|a|\leq q/2$, 
  and let
$f
\in C^2(I;\la_0,\la_1,\la_2,j)$.
Then we have
$$
A_{I}(q;a,b,f) \ll
\frac{1}{|a|}\Big(q+ \la_1  E\Big),
$$
where
\begin{equation}\lab{E}
E=
\frac{|b|^{1/2}q^{1/2}}{\la_2^{1/2}}+
\frac{|b|^{3/2}j\la_2^{1/2}\m(I)}{q^{1/2}}  +
\frac{b^2\la_0}{q}.
\end{equation}
\end{lem}

\begin{proof}
Suppose that $I=[t_1,t_2]$, for $t_1<t_2$.
We begin by following the proof of Lemma
\ref{triv2}.  Thus we may assume that
\eqref{pisa2} holds, with $|1-e(a/q)|^{-1}\ll
|a|^{-1}q$ and $F(t)=\e_q(b f(t))$ for $t_1<t\leq t_2$.
Then it is not hard to conclude that
\begin{equation}\lab{pisa3}
A_{I}(q;a,b,f)
\ll \frac{q}{|a|}\Big(1+|J|\Big),
\end{equation}
where
\begin{align*}
J
&=
\sum_{t_1<n \leq t_2}\e_q(a n)\Big(F(n)-F(n-1)\Big)\\
&=
\sum_{t_1<n \leq t_2}\e_q(a n+b f(n))\Big(1-\e_q(b(f(n-1)-f(n))\Big).
\end{align*}
Let $n \in (t_1,t_2]$.  By the mean value theorem there exists
$\xi \in (n-1,n)$ such that $f(n)-f(n-1)=f'(\xi)$.  Since
$f\in C^2(I;\la_0,\la_1,\la_2,j)$, it follows that
$$
\sup_{t_1<n\leq t_2}|f(n)-f(n-1)| \leq \la_1. %t%
$$
In view of the familiar estimate $\e^{it}=1+it +O(t^2)$, that is
valid for any $t \in \R$, we deduce that
\begin{align*}
1-\e_q(b(f(n-1)-f(n))=& 2 \pi i b\big(f(n)-f(n-1)\big)/q\\
&\quad +
O\big(b^2\la_1|f(n)-f(n-1)|/q^2\big).
\end{align*}
Hence
$$
J\ll \frac{|b|}{q}|S| +
\frac{b^2\la_1}{q^2} \sum_{t_1<n \leq t_2}|f(n)-f(n-1)| \ll
\frac{|b|}{q}|S| +
\frac{b^2\la_0\la_1}{q^2},
$$
where
$$
S=\sum_{t_1<n \leq t_2}\e_q(a n+b f(n))\Big(f(n)-f(n-1)\Big).
$$
Our final task is to handle this sum.

Let $G(t)=f(t)-f(t-1)$ and
$$
T_t=\sum_{t_1<n\leq  t}\e_q(a n + b f(n)),
$$
for any $t\in (t_1,t_2]$.  Then the second derivative
estimate of Van der Corput \cite[Theorem 5.9]{titch} yields
\begin{equation}\lab{p3}
T_t \ll j\m(I)(|b|\la_2/q)^{1/2}+ (|b|\la_2/q)^{-1/2},
\end{equation}
since $f\in C^2(I;\la_0,\la_1,\la_2,j)$ and $t \leq t_2$.
Now an application of partial summation gives
$$
S=T_{t_2}G(t_2)-\int_{t_1}^{t_2} T_tG'(t)\d t.
$$
On applying the mean value theorem to $G$ and $G'$, we therefore
conclude from \eqref{p3} that
$$
S\ll
\Big(j\m(I)(|b|\la_2/q)^{1/2}+ (|b|\la_2/q)^{-1/2}
\Big)\Big(\la_1 + \int_{t_1}^{t_2} |G'(t)| \d t\Big).
$$
But the last integral here is clearly $O(\lambda_1)$, since $f'$ is
monotonic and of constant sign on $(t_1,t_2)$.
Putting all of this together we therefore conclude that \eqref{pisa3}
holds, with
\begin{align*}
qJ \ll
|b||S| +
\frac{b^2\la_0\la_1}{q} \ll \la_1 E,
\end{align*}
in the notation of \eqref{E}.
This completes the proof of the lemma.
\end{proof}

\section{Equidistribution}\lab{equi}

During the course of the proof of Theorem
\ref{main'}, as carried out in \S\S \ref{proof1}--\ref{proof4}
below, we shall need a precise expression for the number of integers
in an interval that lie in a fixed congruence class.
Define the real-valued function $\psi(t)=\{t\}-1/2$, where
$\{t\}$ denotes the fractional part of $t \in \R$.
Then $\psi$ is periodic with period~$1$, and we
have the following simple estimate  \cite[Lemma 3]{math.NT/0412086}.

\begin{lem}\lab{cong}
Let $a, q \in \Z$ be such that $q>0$, and let $t_1, t_2\in \R$ such
that $t_2 \geq t_1$. Then
$$
\#\{t_1<n \leq t_2: n\equiv a\mod{q} \}=\frac{t_2-t_1}{q}+ r(t_1,t_2;a,q),
$$
where
$$
r(t_1,t_2;a,q)=\psi\Big( \frac{t_1-a}{q}\Big) - \psi\Big( \frac{t_2-a}{q}\Big).
$$
\end{lem}

In relation to this result we shall need some control over the average order
of the function $\psi(g(x,y)/q)$, for certain real-valued functions $g$,  as we
range over integers $x,y$ that are restricted to certain intervals
and that satisfy a certain congruence relation modulo $q$.
The simplest scenario is when $g(x,y)$ is actually a polynomial in one
variable, in which case we shall make use of the following result
\cite[Lemma 5]{math.NT/0412086}, established by combining a Fourier
series expansion for $\psi$ with standard bounds for the quadratic
Gauss sum.

\begin{lem}\label{lem:sum_square}%t%
Let $\varepsilon>0$ and let $t \in \R$.  Then for any $a,q \in \Z$ such that
$q>0$ and $\hcf(a,q)=1$, we have
$$
\sum_{\colt{y=1}{\hcf(y,q)=1}}^q
\psi\Big(\frac{t-ay^2}{q}\Big) \ll_\varepsilon q^{1/2+ \varepsilon}.
$$
\end{lem}

We shall also need to examine the average order of $\psi(g(x,y)/q)$
for the more complicated case in which $g(x,y)=f(x)-xy$ for a suitable function $f$.
More precisely, given $a,b,c,q \in \Z$ such that $q>0$ and
$\hcf(abc,q)=1$, and an interval $I \subset \R$, we'll want to study
the sum
\begin{equation}\lab{S_I}
S_I(f,q)=S_I(f,q;a,b,c):=
\sum_{\colt{x\in \Z \cap I}{\hcf(x,q)=1}}\sum_{\colt{y=1}{ay^2\equiv
       bx \mod{q}}}^q \psi\Big(\frac{f(x)-cxy}{q}\Big),
\end{equation}
for suitable real-valued functions $f$ on $I$.  Our estimates for
$S_I(f,q)$ will depend upon the work in the previous section, and we
shall eventually obtain two distinct estimates according to whether we
are in a position to apply Lemma~\ref{triv2} or Lemma~\ref{triv2'}.
We begin however by recording the following ``trivial'' bound for
\eqref{S_I}, which follows from the fact that for fixed integers
$a,b,x$ such that $\hcf(abx,q)=1$, there are
$O_\varepsilon(q^\varepsilon)$ possible
solutions modulo $q$ of the congruence $ay^2\equiv bx (\bmod{~q})$.

\begin{lem}\lab{Lucca}
Let $I \subset \R$ be an interval and
suppose that $a,b,c,q \in \Z$ such that $q>0$ and $\hcf(abc,q)=1$.
Then for any real-valued function $f$ on $I$ we have
$$
S_{I}(f,q)\ll_\varepsilon q^{\varepsilon}\m(I),
$$
where $\m(I)$ is given by \eqref{meas}.
\end{lem}

The starting point for a more sophisticated
treatment of $S_{I}(f,q)$ is the trigonometric formula
\cite{vaaler} for $\psi$, that is due to Vaaler.  For any $t \in \R$,
and any $H \geq 1$,  this implies that
$$
\sum_{0< |h|\leq H}c_h^- {\e}(ht) +O\Big(\frac{1}{H}\Big)
\leq
\psi(t)\leq \sum_{0< |h|\leq H}c_h^+ {\e}(ht) +O\Big(\frac{1}{H}\Big),
$$
for certain coefficients $c_h^-,c_h^+\ll 1/|h|$.  Arguing
as above we therefore deduce that
\begin{equation}\lab{cont3}
S_{I}(f,q)\ll_\varepsilon \frac{q^\varepsilon\m(I)}{H}
+\sum_{h=1}^H\frac{1}{h}\Big|T_{I}(f,q;h)\Big|,
\end{equation}
in the notation of \eqref{meas}, where
$$
T_{I}(f,q;h)
=\sum_{\colt{x \in I\cap\Z}{\hcf(x,q)=1}}\sum_{\colt{y=1}{ay^2\equiv
       bx \mod{q}}}^q \e_q(hf(x)-chxy).
$$
Extending the summation over $x$ to a complete set of residues modulo
$q$, we obtain
\begin{align*}
T_{I}(f,q;h)&=
\sum_{\colt{u=1}{\hcf(u,q)=1}}^q\sum_{x \in I\cap\Z}\frac{1}{q}\sum_{k=1}^q
\e_q(k(u-x))
\hspace{-0.2cm}
\sum_{\colt{v=1}{av^2\equiv
       bu \mod{q}}}^q
\hspace{-0.2cm}
\e_q(hf(x)-chuv)\\
&= \frac{1}{q}\sum_{k=1}^q A_{I}(q;-k,h,f)B(q;h,k),
\end{align*}
where
\begin{align*}
B(q;h,k)&=
\sum_{\colt{u=1}{\hcf(u,q)=1}}^q \sum_{\colt{v=1}{av^2\equiv
       bu \mod{q}}}^q \e_q(ku-chuv)
\end{align*}
and $A_{I}(q;-k,h,f)$ is given by \eqref{exp-sum''}.
By periodicity, we may replace the summation over $1 \leq k \leq q$ by
a summation over $-q/2<k\leq q/2$.

On letting $\overline{b}$ denote the
multiplicative inverse of $b$ modulo $q$, it is easy to see that
$$
B(q;h,k)
= \sum_{\colt{v=1}{\hcf(v,q)=1}}^q \e_q\Big(a\ov{b}(-chv^3+kv^2)\Big).
$$
In order to estimate this sum we must first take
care to remove any possible common factors between $q$ and the
coefficients of $v^3$ and $v^2$. Since $\hcf(abc,q)=1$ by assumption, we
see that
$\hcf(q,a\ov{b}ch,a\ov{b}k)=\hcf(q,h,k),$ whence
$$
T_{I}(f,q;h)= \sum_{d\mid h,q} \frac{1}{dq'}
\sum_{\colt{-q'/2<k'\leq q'/2}{\hcf(k',h',q')=1}} A_{I}(q';-k',h',f)
B(dq';dh',dk').
$$
Here, we have written $h=dh', k=dk'$ and $q=dq'$.

We must now consider the sum $B(dq';dh',dk')$ in more detail.  Each
$v$, modulo $q$, can be written uniquely in the
form $v=y+q'z$ with $1 \leq y \leq q'$ and $1 \leq z \leq d$.
Thus it follows that
\begin{align*}
B(dq';dh',dk')&=\sum_{y=1}^{q'} \sum_{\colt{z=1}{\hcf(y+q'z,dq')=1}}^d
{\e}_{q'}\Big(a\ov{b}(-ch'y^3+k'y^2)\Big)\\
&=\sum_{\colt{y=1}{\hcf(y,q')=1}}^{q'}
{\e}_{q'}\Big(a\ov{b}(-ch'y^3+k'y^2)\Big)
N(d;q',y),
\end{align*}
where $N(d;q',y)$ is the number of positive integers $z\leq d$ for which
$y+q'z$ is coprime to $d$.  But then it is clear that
\begin{align*}
N(d;q',y)
&=\sum_{\ell \mid d}\mu(\ell) \#\{1\leq z \leq d: q'z \equiv -y
\mod{\ell}\}\\
&=\sum_{\colt{\ell \mid d}{\hcf(\ell,q')=1}}\mu(\ell)
\sum_{t=1}^{d/\ell}
\#\{1\leq s \leq \ell: q's \equiv -y \mod{\ell}\}\\
&=d\sum_{\colt{\ell \mid d}{\hcf(\ell,q')=1}}\frac{\mu(\ell)}{\ell}
=f(d,q'),
\end{align*}
say.  In particular we have
$$
f(d,q')= d\phi^*(d)/\phi^*(\hcf(d,q'))\leq d,
$$
where
\begin{equation}\lab{defphi*}
\phi^*(n):=\prod_{p\mid n}(1-1/p) .
\end{equation}
Thus
$
B(dq';dh',dk')=f(d,q')B(q',h',k'),
$
and so
\begin{equation}\lab{dentist}
T_{I}(f,q;h)\ll \sum_{d\mid h,q} \frac{1}{q'}
\sum_{\colt{-q'/2<k'\leq q'/2}{\hcf(k',h',q')=1}} |A_{I}(q';-k',h',f)||B(q';h',k')|.
\end{equation}
We now break the inner sum over $k'$ into two sums: the single term
arising from $k'=0$ and the summation  over $-q'/2<k'\leq q'/2$ such
that $k'\neq 0$.

We begin by handling the overall contribution from the term $k'=0$.
But then it follows from \eqref{triv1} that
$$
A_{I}(q';0,h',f)\ll \m(I),
$$
and from Lemma \ref{exp1} that
$$
B(q';h',0) \ll_\varepsilon {q'}^{2/3+\varepsilon}.
$$
Here we have used the fact that
$\hcf(k',h',q')=\hcf(h',q')=1$.  Combining these two estimates
we therefore obtain the overall contribution
\begin{equation}\lab{cont1}
\ll_\varepsilon q^\varepsilon
\sum_{d\mid h,q} \frac{\m(I)}{(q/d)^{1/3}} \ll_\varepsilon
\frac{q^{2\varepsilon}\hcf(h,q)^{1/3}\m(I)}{q^{1/3}},
\end{equation}
to the right-hand side of \eqref{dentist}.

In order to handle the remaining contribution, our argument bifurcates according to which of
Lemmas \ref{triv2} or \ref{triv2'} we apply to estimate $A_{I}(q';-k',h',f)$.
In either case we may clearly deduce from Lemma \ref{exp3} that
\begin{equation}\lab{app-exp3}
B(q';h',k')\ll_\varepsilon {q'}^{1/2+\varepsilon}\hcf(k',q').
\end{equation}
We begin with an application of Lemma \ref{triv2}, for which we shall
assume that $f \in C^1(I;\la_0)$ for a certain
value of $\la_0\geq 1$.  Thus it follows that
$$
A_{I}(q';-k',h',f) \ll \frac{q'}{k'}\Big(1+\frac{h\la_0}{q}\Big),
$$
since $0<|k'|\leq q'/2$, whence
\begin{align*}
\sum_{\tstack{-q'/2<k'\leq q'/2}{\hcf(k',h',q')=1}{k'\neq 0}}
|A_{I}(q';-k',h',f)||B(q';h',k')|
&\ll_\varepsilon {q'}^{3/2+2\varepsilon}\Big(1+\frac{h\la_0}{q}\Big).
\end{align*}
Here we have used the trivial observation that
\begin{equation}\lab{jung}
\sum_{1\leq a\leq A}\frac{\hcf(a,b)}{a} \leq \sum_{d \mid b} d
    \sum_{1\leq a \leq
    A/d}\frac{1}{ad} \ll \tau(b)\log A,
\end{equation}
for any $A\geq 2$ 
and any $b \in \N$, together with the upper bound
$\tau(n)=O_\varepsilon (n^{\varepsilon})$ for the divisor function.
We therefore obtain the overall contribution
\begin{equation}\lab{cont2}
\ll_\varepsilon \sum_{d\mid h,q}
q^{1/2+2\varepsilon}\Big(1+\frac{h\la_0}{q}\Big)
\ll_\varepsilon q^{1/2+3\varepsilon}\Big(1+\frac{h\la_0}{q}\Big),
\end{equation}
to the right-hand side of \eqref{dentist} from this case.
Alternatively, let $j,\la_0,\la_1,\la_2 \in \R$ such that \eqref{la_}
holds, and suppose that $f \in
C^2(I;\la_0,\la_1,\la_2,j)$.  Then it follows that
$$
A_{I}(q';-k',h',f) \ll
\frac{q'}{k'}(1+\la_1 E),
$$
where
\begin{equation}\lab{E'}
E=
\frac{h^{1/2}}{\la_2^{1/2}q^{1/2}}+
\frac{h^{3/2}j\la_2^{1/2}\m(I)}{{q}^{3/2}}  +
\frac{{h}^2\la_0 }{{q}^2}.
\end{equation}
We may combine this with \eqref{app-exp3} and  \eqref{jung}
to obtain the overall contribution
\begin{equation}\lab{cont2'}
\sum_{d\mid h,q} \frac{1}{q'}
\sum_{\tstack{-q'/2<k'\leq q'/2}{\hcf(k',h',q')=1}{k'\neq 0}}
|A_{I}(q';-k',h',f)||B(q';h',k')|
\ll_\varepsilon {q}^{1/2+3\varepsilon}(1+\la_1 E)
\end{equation}
to the right-hand side of \eqref{dentist} from this case.

Let us begin by drawing together
\eqref{cont1} and \eqref{cont2} in
\eqref{dentist}, before
then inserting the resulting bound into \eqref{cont3}.  In view of
\eqref{jung}  we have shown that
\begin{align*}
S_{I}(f,q)
&\ll_\varepsilon
\frac{q^\varepsilon\m(I)}{H}
+q^{3\varepsilon}\sum_{h=1}^H\Big(
\frac{\hcf(h,q)\m(I)}{hq^{1/3}}
+\frac{q^{1/2}}{h}+ \frac{\la_0}{q^{1/2}}\Big)\\
&\ll_\varepsilon
q^{3\varepsilon}H^{\varepsilon}\Big(
\frac{\m(I)}{H}
+\frac{\m(I)}{q^{1/3}}
+q^{1/2}+H\frac{\la_0}{q^{1/2}} \Big),
\end{align*}
for any $f \in C^1(I;\la_0)$  and any $H \geq 1$.
Suppose first that $\m(I)q^{1/2} \geq \la_0$.  Then we may select
$$
H=\frac{\m(I)^{1/2}q^{1/4}}{\la_0^{1/2}},
$$
to get
$$
S_{I}(f,q)
\ll_\varepsilon
q^{4\varepsilon}\m(I)^\varepsilon\Big(
\frac{\m(I)}{q^{1/3}}
+q^{1/2}+
\frac{\m(I)^{1/2}\la_0^{1/2}}{q^{1/4}}\Big).
$$
Alternatively, if $\m(I)q^{1/2} \leq \la_0$ we employ the trivial
estimate Lemma \ref{Lucca} for $S_{I}(f,q)$, to conclude that
$$
S_{I}(f,q) \ll_\varepsilon
q^{\varepsilon}\m(I) \ll_\varepsilon
q^\varepsilon
\frac{\m(I)^{1/2}\la_0^{1/2}}{q^{1/4}}.
$$
On combining these two estimates and redefining the choice of $\varepsilon$,
we have therefore established the following result.

\begin{lem}\lab{cor1}
Let $I \subset \R$ be an  interval and let $\la_0 \geq 1$.
Suppose that $a,b,c,q \in \Z$ such that $q>0$ and $\hcf(abc,q)=1$,
and let $f \in C^1(I;\la_0)$.  Then we have
$$
S_{I}(f,q)\ll_\varepsilon
q^{\varepsilon}\m(I)^\varepsilon\Big(
q^{1/2}+\frac{\m(I)}{q^{1/3}}
+ \frac{\la_0^{1/2}\m(I)^{1/2}}{q^{1/4}}\Big),
$$
where $\m(I)$ is given by \eqref{meas}.
\end{lem}

We may obtain an alternative estimate for $S_I(f,q)$ by drawing
together  \eqref{cont1} and \eqref{cont2'} in
\eqref{dentist}, when
$f \in C^2(I;\la_0,\la_1,\la_2,j)$ for $j,\la_0,
\la_1,\la_2\in \R$ such that~\eqref{la_}
holds. On inserting the resulting estimate for \eqref{dentist} into
\eqref{cont3} we conclude that
\begin{align*}
S_{I}(f,q)
&\ll_\varepsilon
\frac{q^\varepsilon\m(I)}{H}
+q^{3\varepsilon}\sum_{h=1}^H\Big(
\frac{\hcf(h,q)\m(I)}{hq^{1/3}}
+\frac{q^{1/2}}{h}+ \frac{q^{1/2}\la_1}{h}  E\Big),
\end{align*}
where $E$ is given by \eqref{E'}.
But then \eqref{jung} yields
$$
S_{I}(f,q)
\ll_\varepsilon
q^{3\varepsilon}H^{\varepsilon}\Big(
\frac{\m(I)}{H}
+\frac{\m(I)}{q^{1/3}}
+q^{1/2}+ F\Big),
$$
where
$$
F=
\frac{H^{1/2}\la_1}{\la_2^{1/2}}+
\frac{H^{3/2}j\la_1\la_2^{1/2}\m(I)}{q}+
\frac{H^2\la_0\la_1}{q^{3/2}}.
$$
Suppose first that $\la_2\m(I)^2 \geq \la_1^2$. Then we may select
$$
H=\frac{\la_2^{1/3}\m(I)^{2/3}}{\la_1^{2/3}},
$$
and it follows that
\begin{align*}
S_I(f,q)
&\ll_\varepsilon
q^{4\varepsilon}\Big(
q^{1/2}+ \frac{\m(I)}{q^{1/3}}+
\frac{\la_1^{2/3}\m(I)^{1/3}}{\la_2^{1/3}}+
\frac{j\la_2\m(I)^2}{q}+
\frac{\la_0 \la_2^{2/3}\m(I)^{4/3}}{\la_1^{1/3}q^{3/2}}
\Big).
\end{align*}
Alternatively, if $\la_2\m(I)^2 \leq \la_1^2$ then Lemma \ref{Lucca}
implies that
$$
S_{I}(f,q) \ll_\varepsilon
q^{\varepsilon}\m(I) \ll_\varepsilon
q^\varepsilon
\frac{\la_1^{2/3}\m(I)^{1/3}}{\la_2^{1/3}}.
$$
On combining these two estimates and redefining the choice of $\varepsilon$,
we have therefore established the following result.

\begin{lem}\lab{cor2}
Let $I \subset \R$ be an  interval and let $j,\la_0,\la_1,\la_2 \in
\R$ such that \eqref{la_} holds.
Suppose that $a,b,c,q \in \Z$ such that $q>0$ and $\hcf(abc,q)=1$, and
let $f \in C^2(I;\la_0,\la_1,\la_2,j)$. Then we have
\begin{align*}
S_I(f,q)&\ll_\varepsilon
q^{\varepsilon}\Big(
q^{1/2}+ \frac{\m(I)}{q^{1/3}}+
\frac{\la_1^{2/3}\m(I)^{1/3}}{\la_2^{1/3}}+
\frac{j\la_2\m(I)^2}{q}+
\frac{\la_0 \la_2^{2/3}\m(I)^{4/3}}{\la_1^{1/3}q^{3/2}}
\Big),
\end{align*}
where $\m(I)$ is given by \eqref{meas}.
\end{lem}

\section{The real-valued functions $g_1$ and $g_2$}\lab{real}

The purpose of this section is to introduce and
analyse a number of real-valued functions that play a
pivotal role in subsequent sections. In fact they will arise in \S
\ref{proof2} as boundary curves for the heights of the variables to be
introduced during our passage to the universal torsor. It is
precisely to some of these functions that we will
ultimately apply the results of the previous
section.

We begin by introducing a function $g_1:[0,1] \rightarrow \R$ on the
unit interval, given by
\begin{equation}\lab{g1}
g_1(v):=-(\min\{ 1/v^{4}, 1+1/v^2\})^{1/3}.
\end{equation}
Next we introduce functions
$g_{21},g_{22}:(-\infty,1]\times [0,1] \rightarrow \R$, which are
given  by
\begin{equation}\lab{g21}
g_{21}(u,v):=
\left\{
\begin{array}{ll}
0, & \mbox{if $-1\leq u \leq 1$},\\
\sqrt{-1-u^3}, & \mbox{if $u \leq -1$},
\end{array}
\right.
\end{equation}
and
\begin{equation}\lab{g22}
g_{22}(u,v):=
\left\{
\begin{array}{ll}
\sqrt{1-u^3}, & \mbox{if $-(1/v^2-1)^{1/3}\leq u \leq 1$},\\
1/v, & \mbox{if $u \leq -(1/v^2-1)^{1/3}$},
\end{array}
\right.
\end{equation}
respectively.  Finally let us define the function
$g_{2}:\R^2 \rightarrow \R$, by
\begin{equation}\lab{g2'}
g_{2}(u,v):=
\left\{
\begin{array}{ll}
g_{22}(u,v)-g_{21}(u,v),
& \mbox{if $g_1(v)\leq u\leq 1$ and $v \in [0,1]$,}\\
0, & \mbox{otherwise}.
\end{array}
\right.
\end{equation}
Then we clearly have
\begin{equation}\lab{eq:integral_gone}
\begin{split}
g_2(u,v)&=\sqrt{\min \{ 1/v^2, 1-u^3\}}-\sqrt{\max \{ 0, -1-u^3\}}\\
&=
\int_{\{t \in \R : ~0 \leq tv \le 1, ~|t^2+u^3| \le 1\}} \dd t,
\end{split}
\end{equation}
for $g_1(v)\leq u\leq 1$ and $v \in [0,1]$.

Write $D_i g_2(u,v)$ for the partial derivative of $g_2$ with respect
to the $i$-th variable, for $i=1,2$.
Then we shall need to know something about the behaviour of
$D_1 g_2(u,v)$ and $|D_1 g_2(u,v)|$ as we integrate over all values of
$u$. For this it will be necessary to calculate $g_2(u,v)$
explicitly.  It suffices to restrict attention to $u,v$ such that
$g_1(v)\leq u\leq 1$ and $v \in [0,1]$, since $g_2(u,v)$ is defined to
be zero in all other cases.
Thus one easily combines \eqref{g21}, \eqref{g22} and
\eqref{eq:integral_gone} to deduce that
$$
g_2(u,v)=\sqrt{ 1-u^3},
$$
if $0\leq u\leq 1$.  If $2^{-1/2}<v\leq 1$ and $1/v^2-1<-u^3\leq 1$, then
$$
g_2(u,v)=1/v.
$$
If $0\leq  -u^3\leq \min\{1/v^2-1,1\}$, we have
$$
g_2(u,v)=\sqrt{ 1-u^3},
$$
while if $0<v\leq 2^{-1/2}$  and $1\leq -u^3\leq   1/v^2-1 $, we have
$$
g_2(u,v)=\sqrt{ 1-u^3}-\sqrt{ -1-u^3}=\frac{2}{\sqrt{ 1-u^3}+\sqrt{
-1-u^3} }.
$$
Finally,
$$
g_2(u,v)=1/v -\sqrt{ -1-u^3}
$$
if $\max\{1/v^2-1,1\}<-u^3\leq \min\{1/v^2+1,1/v^{4}\}$.
It follows from these calculations that
\begin{equation}\lab{intD1g1-a}
\int_{-\infty}^\infty  D_1g_2(u,v) \d u=g_2(1,v)-g_2(g_1(v),v)=-g_2(g_1(v),v),
\end{equation}
since $g_2(1,v)=0$ for any $v\in \R$.  Furthermore it is now
straightforward to conclude that
\begin{equation}\lab{intD1g1-b}
\int_{-\infty}^\infty  |D_1g_2(u,v)| \d u \ll 1,
\end{equation}
for any $v \in \R$.
The final fact that we shall need to highlight here is the elementary equality
\begin{equation}\lab{g=0}
g_2(g_1(v),v)=0, \qquad (0<v\leq 2^{-1/2}).
\end{equation}

\section{Proof of Theorem \ref{main'}: preliminaries}\lab{proof1}

We are now ready to commence the proof of Theorem \ref{main'} in
earnest, for which it is necessary to introduce some more notation.
For any $n \geq 2$ we let $Z^{n+1}$ denote the set of
primitive vectors in $\Z^{n+1}$,
where $\ma{v}=(v_0,\ldots,v_n) \in \Z^{n+1}$ is
said to be primitive if $\hcf(v_0,\ldots,v_n)=1.$
Finally we shall let $\Z_*^{n+1}$ (resp. $Z_*^{n+1}$) denote the
set of vectors $\ma{v} \in \Z^{n+1}$ (resp. $\ma{v} \in Z^{n+1}$)
such that $v_0\cdots v_n \neq 0$.

In this section we shall establish a preliminary estimate for
$\NU(B)$, which paves the way towards the universal torsor
calculation in the following section.
If $x=[\x]\in \bfP^3(\Q)$ is represented by
the vector $\x \in Z^4$, then $H(x)=|\x|=\max_{0\leq
i\leq 3 }|x_i|$.  Thus it is easy to see that
$$
N_{U,H}(B)= \frac{1}{2}\#\{\x\in Z^4: |\x| \leq
B, ~x_1x_2^2+x_2x_0^2+x_3^3=0\},
$$
since $\x$ and $-\x$ represent the same point in $\bfP^3$.
We proceed by considering the contribution to $\NU(B)$ from vectors
$\x \in Z^4$ which contain zero components.
Define the set
\begin{equation}\lab{MB}
\E(B):=\{\x\in Z_*^4: ~x_0,x_2>0,~ |\x| \leq B, ~f(\x)=0\},
\end{equation}
where $f(\x)=x_1x_2^2+x_2x_0^2+x_3^3$.
Then we have the following result.

\begin{lem}\lab{Reduc1}
Let $B\geq 1$.  Then we have
$$
N_{U,H}(B)=2\#\E(B) + \frac{12}{\pi^2}B+O(B^{2/3}).
$$
\end{lem}
\begin{proof}
Suppose that $\x\in Z^4$ is a vector such that
$$
x_0x_1=0, \quad |x_0|,|x_1|,|x_2|,|x_3|\leq B,
$$
and $x=[x_0,\ldots,x_3]\in S$, where $S$ denotes
the surface $f(\x)=0$. If $x_0=0$ then
$x_1x_2^2+x_3^3=0$, and so we are interested in vectors of the shape
$(x_1,x_2,x_3)=\pm (a^3, b^3, -ab^2)$ for coprime $a,b \in \N$.
Such points therefore contribute $O(B^{2/3})$ to $\NU(B)$.
The case in which $x_1=0$ is similar.  Suppose now that $x \in S$ is
represented by a vector $\x\in Z^4$ such that
$$
x_2=0, \quad |x_0|,|x_1|,|x_3|\leq B.
$$
Then necessarily $x_3=0$ and the point is not to be counted, since it
lies on the unique line in $S$.  Finally
we suppose that $\x\in Z^4$ satisfies
$$
x_3=0, \quad 0< |x_0|,|x_1|,|x_2|\leq B,
$$
and $[x_0,x_1,x_2,0]\in S$.  Then we must have $x_0^2+x_1x_2=0$,
so that $\x=(\pm ab,-a^2,b^2,0)$
or $\x=(\pm ab,a^2,-b^2,0)$ for coprime $a, b \in \N$.
Now the number of coprime $a,b\in \N$ such that $a,b\leq Y$ is
$6Y^2/\pi^2 +O(Y\log Y)$. %t%
Hence the overall contribution from this case is
$12B/\pi^2+O(B^{1/2})$, whence
$$
N_{U,H}(B)=\frac{1}{2}
\#\{\x\in Z_*^4: |\x| \leq B, ~f(\x)=0\}
+ \frac{12}{\pi^2}B+O(B^{2/3}).
$$
We complete the proof of Lemma \ref{Reduc1} by choosing $x_0>0$ and
$x_2>0$, as we clearly may.
\end{proof}

\section{Proof of Theorem \ref{main'}: the universal torsor}\label{proof2}

The purpose of this section is to establish a bijection between the
rational points on the open subset $U$ of the cubic surface $S$, and the
integral points on the universal torsor
above $\tS$, which are subject to a number of coprimality conditions.
Along the way we shall introduce new variables
$$
\xi_1,\xi_2,\xi_3,\xi_\ell,\xi_4,\xi_5,\xi_6,\tau_1,\tau_2,\tau_\ell,
$$
and it will be convenient to henceforth write
$\bxi=(\xi_1,\xi_2,\xi_3,\xi_\ell,\xi_4,\xi_5,\xi_6)$ and
$\bt=(\tau_1,\tau_2,\tau_\ell)$.
Furthermore we shall make frequent use of the notation
$$
\base {n_1}{n_2}{n_3}{n_\ell}{n_4}{n_5}{n_6} :=
\xi_1^{n_1}\xi_2^{n_2}\xi_3^{n_3} \xi_\ell^{n_\ell}
\xi_4^{n_4}\xi_5^{n_5}\xi_6^{n_6},
$$
in all that follows.

Hassett and Tschinkel \cite[\S 3]{MR2029868}  have  calculated the Cox
ring of $\tS$ as being given by
$$
\rom{Cox}(\tilde{S})=\Spec(\Q[\bxi,\bt]/(T(\bxi,\bt))),
$$
where
$$
T(\bxi,\bt) :=
\tau_\ell\xi_\ell^3\xi_4^2\xi_5+\tau_2^2\xi_2+\tau_1^3\xi_1^2\xi_3.
$$
The universal torsor $\TT$ is an open subset of the affine hypersurface
\begin{equation}\lab{ut}
\tau_\ell\xi_\ell^3\xi_4^2\xi_5+\tau_2^2\xi_2+\tau_1^3\xi_1^2\xi_3=0,
\end{equation}
together with a map $\Psi: \TT \to S$,
given by
\begin{equation}\label{eq:sub}
\begin{split}
      \Psi^*(x_0) &= \base 1 2 2 0 1 2 3 \tau_2\\
      \Psi^*(x_1) &= \tau_\ell\\
      \Psi^*(x_2) &= \base 2 3 4 3 4 5 6\\
      \Psi^*(x_3) &= \base 2 2 3 1 2 3 4 \tau_1.
\end{split}
\end{equation}
Below we shall use our own methods to show how $\Psi$ gives a bijection
between $U(\Q)$ and the integral points on the hypersurface \eqref{ut}, subject to certain
coprimality conditions.
Recall the notation introduced in \S \ref{sec:manins_conjecture} to describe the
Picard group $\Pic(\tS)$ of $\tS$.  Then it is interesting to remark that the
variables $\xi_1,\ldots,\xi_6,\xi_\ell$ arise as non-zero
sections that generate $\Gamma(E_1),\ldots,\Gamma(E_6),\Gamma(\ell)$,
respectively,
where $\Gamma(A)=H^0(\tS,A)$ for any divisor $A\in \Pic(\tS)$.
The variables $\tau_1,\tau_2, \tau_\ell$ are certain extra sections
that are needed to generate the full Cox ring.
These are chosen to be sections of
the nef divisor class that is dual
to $F_1, F_2, \ell$, respectively.  The interested reader should
consult the calculation of Hassett and Tschinkel for further details.

We now turn to the task of establishing the bijection.
Let us define the set
\begin{equation}\lab{M-infty}
\E:=\{\x\in Z_*^4: x_0,x_2>0, ~x_1x_2^2+x_2x_0^2+x_3^3=0\},
\end{equation}
so that $\E=\E(\infty)$ in the notation of \eqref{MB}.
Then we begin by demonstrating a bijection between this set and the set
$\TT_1$ of $(\bxi,\bt)\in \N^7\times \Z_*\times \N \times \Z_*$
satisfying $T(\bxi,\bt)=0$, such that
\begin{equation}\lab{co-1}
|\mu(\xi_1\xi_2\xi_3\xi_4\xi_5)|=1,
\quad \hcf(\tau_1,\xi_2\xi_\ell\xi_4\xi_5)=1,
\end{equation}
and
\begin{equation}\lab{co-2}
\hcf(\tau_2,\xi_1\xi_3)=\hcf(\tau_\ell,\xi_4\xi_5\xi_6)=1.
\end{equation}
This is achieved in the following result.

\begin{lem}\label{lem:bijection_different_coprime}
The map $\Psi$ induces a bijection between $\TT_1$ and $\E$.
\end{lem}

\begin{proof}
The idea of the proof is simply to go through a series of elementary
considerations about the divisibility properties
enjoyed by elements of the set~$\E$.  At each
stage we shall replace the original variables by products
of new ones which fulfil certain
auxiliary conditions, and which will be uniquely determined by the process.

To begin with, let $\x \in \E$ and note that $x_2 \mid x_3^3$. We may
therefore write $x_2=y_1y_2^2y_3^3$ with $y_1,y_2,y_3 \in \N$ such
that $|\mu(y_1y_2)|=1$, where each triple occurrence of a prime factor
of $x_2$ is put in $y_3$ and each double occurrence in~$y_2$. 
Then $x_3 = y_1y_2y_3z$ must hold for a suitable $z \in \Z_*$.
Substituting into $f$ and dividing by $y_1y_2^2y_3^3$ gives the new equation
\[
f_1(x_0,x_1,y_1,y_2,y_3,z) = x_1y_1y_2^2y_3^3 + x_0^2 + y_1^2y_2z^3 =0.
\]
Now $y_1y_2\mid x_0^2$, and since $|\mu(y_1y_2)|=1$, we have
$y_1y_2\mid x_0$.
Write $x_0=y_1y_2w$ for a suitable $w\in \N$. Substituting, and dividing by
$y_1y_2$, we therefore obtain
\[
f_2(x_1,y_1,y_2,y_3,z,w) = x_1y_2y_3^3+w^2y_1y_2+y_1z^3 = 0.
\]
Since $y_2\mid y_1z^3$ and $|\mu(y_1y_2)|=1$, we must have $y_2\mid z$. Writing
$z=y_2z'$, where $z' \in \Z_*$, yields
\[
f_3(x_1,y_1,y_2,y_3,w,{z'}) = x_1y_3^3+w^2y_1+y_1y_2^2{z'}^3 = 0,
\]
after dividing by $y_2$.
Since $y_1$ divides our original variables $x_0,x_2,x_3$, it cannot
divide $x_1$. Once combined with the fact that $y_1$ is square-free,
the fact $y_1\mid x_1y_3^3$ implies
that $y_1\mid y_3$. Hence $y_3 = y_1y_3'$, where $y_3' \in \N$, and we obtain
\[
f_4(x_1,y_1,y_2,w,{z'},{y_3'}) = x_1y_1^2{y_3'}^3+w^2+y_2^2{z'}^3 = 0.
\]
Let $a=\hcf({y_3'},{z'}) \in \N$ and write ${y_3'}=a{y_3''}$ and ${z'}=a{z''}$,
for ${y_3''}\in \N$ and ${z''} \in \Z_*$. This gives
\[
f_5(x_1,y_1,y_2,w,{z''},{y_3''},a) =
x_1y_1^2{y_3''}^3a^3+w^2+y_2^2{z''}^3a^3= 0.
\]
Now $a^3\mid w^2$. Writing $a=\xi_6^2\xi_2$, for $\xi_2,\xi_6 \in \N$
such that $|\mu(\xi_2)|=1$, gives $w=w'\xi_6^3\xi_2^2$ for suitable
$w'\in \N$.  This leads to the equation
\[
f_6(x_1,y_1,y_2,{z''},{y_3''},{w'},\xi_2,\xi_6) =
    x_1y_1^2{y_3''}^3+{w'}^2\xi_2+y_2^2{z''}^3=0.
\]
Let $\xi_5 = \hcf({y_3''},{w'}) \in \N$ and write ${y_3''}={\xi_\ell}\xi_5$ and
${w'}=w''\xi_5$, for suitable ${\xi_\ell},w'' \in \N$. Then
\[
f_7(x_1,y_1,y_2,{z''},w'',\xi_2,{\xi_\ell},\xi_5,\xi_6) =
x_1y_1^2{\xi_\ell}^3\xi_5^3 + {w''}^2\xi_2\xi_5^2 + y_2^2{z''}^3 = 0.
\]
Since $\hcf({y_3''},{z''}) = 1$, we also have
$\hcf({\xi_\ell}\xi_5,{z''}) = 1$.
Therefore $\xi_5^2\mid y_2^2$, and so $\xi_5\mid y_2$.  We proceed to write
$y_2=\xi_1\xi_5$, with $\xi_1\in \N$, and so obtain
\[
f_8(x_1,y_1,{z''},w'',\xi_1,\xi_2,{\xi_\ell},\xi_5,\xi_6) =
x_1y_1^2{\xi_\ell}^3\xi_5+{w''}^2\xi_2+\xi_1^2{z''}^3 = 0.
\]
Let $\xi_3 = \hcf(w'',y_1) \in \N$. Since
$|\mu(y_1y_2)|=1$, we have $\hcf(\xi_1,\xi_3) = 1$. Therefore,
$\xi_3 \mid {z''}^3$ and even $\xi_3\mid {z''}$.  Write $w''=\tau_2 \xi_3$
for suitable $\tau_2 \in \N$, and $y_1=\xi_4\xi_3$ for suitable $\xi_4
\in \N$, and finally ${z''}=\tau_1 \xi_3$ for suitable $\tau_1 \in
\Z_*$. 
Then on writing $x_1=\tau_\ell \in \Z_*$, we therefore obtain
the torsor equation $T(\bxi,\bt)=0$, as given by \eqref{ut}.

Now it is easy to check that the substitutions we have made lead to a
bijection between the set of
$\x\in \Z_*^4$ such that
$$
x_0,x_2>0, \quad x_1x_2^2+x_2x_0^2+x_3^3=0,
$$
and the set  of
$(\bxi,\bt)\in \N^7\times \Z_*\times \N \times \Z_*$
satisfying $T(\bxi,\bt)=0$,
and that this bijection is given  by \eqref{eq:sub}.
In order to complete the proof of the lemma it therefore remains to
collect together the coprimality relations that are satisfied
by the new variables $\bxi, \bt$,  and consider how they correspond to
the primitivity of the vectors in the set $\E$ under the map $\Psi$.

But on tracing through our argument, one easily checks that we have
imposed the coprimality conditions
\begin{align*}
&|\mu(y_1y_2)|=|\mu(\xi_1\xi_3\xi_4\xi_5)|=1,\qquad |\mu(\xi_2)|=1,
\qquad \hcf(\tau_2,\xi_4)=1,\\
      &\hcf({y_3''},{z''}) = \hcf({\xi_\ell}\xi_5,\tau_1\xi_3) = 1,\qquad
      \hcf({\xi_\ell},w'')
      = \hcf({\xi_\ell},\tau_2\xi_3) = 1,
\end{align*}
on $\bxi, \bt$, together with the condition
$\hcf(\tau_\ell,\xi_1\xi_2\xi_3\xi_4\xi_5\xi_6)=1$ that arises from the
condition $\hcf(x_0,x_1,x_2,x_3)=1$. A little thought reveals that
these conditions can be rewritten as
\begin{align*}
&|\mu(\xi_1\xi_3\xi_4\xi_5)|=|\mu(\xi_2)|=\hcf(\xi_3,\xi_\ell)=1,\\
&\hcf(\tau_1,\xi_\ell\xi_5)=\hcf(\tau_2,\xi_\ell\xi_4)=\hcf(\tau_\ell,\xi_1\xi_2\xi_3\xi_4\xi_5\xi_6)=1.
\end{align*}
Our final task is therefore to show, that once taken together with the
equation \eqref{ut}, these coprimality relations are equivalent
to the conditions \eqref{co-1} and \eqref{co-2} that
appear in the definition of $\TT_1$.
We content ourselves with checking the reverse implication, the
other direction being entirely similar.  For this it is clearly enough to show
that the relations \eqref{ut}, \eqref{co-1} and \eqref{co-2}
together imply that
$$
\hcf(\xi_3,\xi_\ell)=
\hcf(\tau_2,\xi_\ell\xi_4)=\hcf(\tau_\ell,\xi_1\xi_2\xi_3)=1.
$$
Suppose first that there is a common prime factor $p$ of
$\xi_3,\xi_\ell$ in \eqref{ut}.  Then $p \mid \tau_2^2\xi_2$, which is
impossible by \eqref{co-1} and \eqref{co-2}.  Hence
$\hcf(\xi_3,\xi_\ell)=1$.  If there is a prime $p$ dividing
$\tau_2$ and $\xi_\ell\xi_4$ in \eqref{ut}, then $p$ divides
$\tau_1^3\xi_1^2\xi_3$, which is also impossible by \eqref{co-1} and
\eqref{co-2}. Finally one deduces that $\tau_\ell$ must be coprime to
$\xi_1\xi_2\xi_3$ in a similar fashion.
This completes the proof of Lemma~\ref{lem:bijection_different_coprime}.
\end{proof}

During the course of our work it will be convenient to replace
the coprimality relations \eqref{co-1} in the definition of $\TT_1$ by
relations of the shape
\begin{equation}\lab{eq:cricket}
|\mu(\xi_2\xi_3\xi_4\xi_5)|=1,
\quad \hcf(\xi_1,\xi_2)=1,
\end{equation}
and 
\begin{equation}\lab{co-1'}
\hcf(\tau_1,\xi_2\xi_3\xi_\ell\xi_4\xi_5\xi_6)=1.
\end{equation}
Thus let us define $\TT_2$ to be the set of $(\bxi,\bt)\in \N^7\times
\Z_*\times \N \times \Z_*$ satisfying $T(\bxi,\bt)=0$,
such that \eqref{co-2}, \eqref{eq:cricket} and \eqref{co-1'}
holds.  On recalling the definition (\ref{M-infty}) of $\E$, we have
the following result.

\begin{lem}\label{bij}
The map $\Psi$ induces a bijection between $\TT_2$ and $\E$.
\end{lem}

\begin{proof}
In view of Lemma \ref{lem:bijection_different_coprime} it will suffice to
establish a suitable bijection between $\TT_1$ and $\TT_2$.
Let $(\bxi,\bt)\in \TT_1$.  We decompose the
coordinates into their prime factors by writing
$$
\xi_i=\prod_p p^{m_{i,p}}, \quad \tau_1=\pm \prod_{p}p^{n_{1,p}},
\quad \tau_2=\prod_{p}p^{n_{2,p}},
\quad \tau_\ell=\pm \prod_{p}p^{n_{\ell,p}}
$$
for $i\in\{1,2,3,\ell,4,5,6\}$.
We shall
henceforth set $m_i=m_{i,p} $ and $n_j=n_{j,p}$, for convenience.
It clearly follows from \eqref{co-1} and \eqref{co-2} that
\begin{equation}\lab{A}
\begin{array}{c}
m_1+m_2+m_3+m_4+m_5 \in \{0,1\},\\
\min\{n_1,m_2+m_\ell+m_4+m_5\}=0,\\
\min\{n_2,m_1+m_3\}=\min\{n_\ell,m_4+m_5+m_6\}=0,
\end{array}
\end{equation}
and once combined with \eqref{ut}, it is easy to deduce the further relation
\begin{equation}\lab{A'}
\min\{m_6,n_1,n_2\}=0.
\end{equation}
We now construct a map $\Phi:\TT_1\rightarrow \TT_2$.  This is obtained
via $\Phi: (\bxi,\bt) \mapsto
(\xi_1',\ldots,\xi_6',\tau_1',\tau_2',\tau_\ell')=(\bxi',\bt')$, where
$$
\xi_i'=\prod_p p^{m_{i,p}'}, \quad \tau_1'=\pm \prod_{p}p^{n_{1,p}'},
\quad \tau_2'=\prod_{p}p^{n_{2,p}'},
\quad \tau_\ell'=\pm \prod_{p}p^{n_{\ell,p}'}
$$
for suitable values of $m_i'=m_{i,p}'$ and $n_j'=n_{j,p}'.$
The idea is to transfer powers of primes $p$ between $\xi_6$ and the
variables appearing in the final monomial $\tau_1^3\xi_1^2\xi_3$ of \eqref{ut}.
We shall set
$m_{i}'=m_{i}$ for $i \in \{2,\ell,4,5\}$, and
$n_{j}'=n_{j}$ for $j \in \{2,\ell\}$.  The remaining values of
$m_{i}',n_{j}'$ are determined in the following way, in which $k$ denotes
an arbitrary non-negative integer.
\begin{itemize}
\item
     If $n_{1} = 2k+1, m_{6} \geq k+1, m_{3} = 0$, then
     \[(m_{1}',m_{3}',m_{6}',n_{1}') = (m_{1}+3k+1,1,m_{6}-k-1,n_{1}-2k-1).\]
\item
     If $n_{1} = 2k+1, m_{6} \geq k+1, m_{3} = 1$, or
     if $n_1>2k, m_{6} =k, m_{3} = 1$, then
     \[(m_{1}',m_{3}',m_{6}',n_{1}') =
     (m_{1}+3k+2,0,m_{6}-k,n_{1}-2k-1).\]
\item
If $n_{1} > 2k, m_{6} =k, m_{3} = 0$, or
if $n_{1} = 2k, m_{6} \geq k,$ then
     \[(m_{1}',m_{3}',m_{6}',n_{1}') =
     (m_{1}+3k,m_{3},m_{6}-k,n_{1}-2k).\]
\end{itemize}
We note that the five possibilities detailed here are exhaustive.
Furthermore, on recalling the conditions \eqref{A} and \eqref{A'}, it
follows from these definitions that
$$
\begin{array}{c}
m_2'+m_3'+m_4'+m_5' \in \{0,1\}, \quad \min\{m_1',m_2'\}=0,\\
\min\{n_1',m_2'+m_3'+m_\ell'+m_4'+m_5'+m_6'\}=0,\\
\min\{n_2',m_1'+m_3'\}=\min\{n_\ell',m_4'+m_5'+m_6'\}=0,
\end{array}
$$
and so $\bxi',\bt'$ satisfy the coprimality relations \eqref{co-2}, \eqref{eq:cricket}
and \eqref{co-1'}.  Finally, it is easily checked that $(\bxi',\bt')\in \TT_2$.

Conversely, let $(\bxi',\bt')\in \TT_2$.
Then we construct a map $\Phi^{-1}:\TT_2\rightarrow \TT_1$, given by
$\Phi^{-1}: (\bxi',\bt') \mapsto (\bxi,\bt)$, by setting
$m_{i}=m_{i}'$ for $i \in \{2,\ell,4,5\}$, and
$n_{j}=n_{j}'$ for $j \in \{2,\ell\}$.  The remaining values of
$m_{i},n_{j}$ depend on the value of $m_{1}'$ modulo $3$, and are
determined in the following way, in which $k$ again denotes an
arbitrary non-negative integer.
\begin{itemize}
\item
If $m_{1}' \in \{3k+1, 3k+2\}$ and $m_{3}' = 1$, then
\[(m_{1},m_{3},m_{6},n_{1}) =
     (m_{1}'-3k-1,0,m_{6'}+k+1,n_{1}'+2k+1).\]
\item
If $m_{1}' = 3k+2$ and $m_{3}' = 0$, then
     \[(m_{1},m_{3},m_{6},n_{1}) =
     (m_{1}'-3k-2,1,m_{6}'+k,n_{1}'+2k+1).\]
\item
If $m_{1}' = 3k+1$ and $m_{3}' = 0$, or if $m_{1}'=3k$, then
     \[(m_{1},m_{3},m_{6},n_{1}) =
     (m_{1}'-3k,m_{3}',m_{6}'+k,n_{1}'+2k).\]
\end{itemize}
Again it is a straightforward to check that
$(\bxi,\bt)\in \TT_1$.

Finally one verifies that $\Phi\Phi^{-1}=\rom{id}_{\TT_2}$ and $
\Phi^{-1}\Phi=\rom{id}_{\TT_1}$,
and that the points $(\bxi,\bt)\in \TT_1$ and $(\bxi',\bt')\in \TT_2$
correspond to the same point in $\E$ under the map $\Psi$.
This therefore completes the proof of Lemma \ref{bij}.
\end{proof}

We are now ready to relate the quantity $\#\E(B)$, as given by (\ref{MB}),
to an appropriate count on the universal torsor.
Now it is clear that under the substitutions (\ref{eq:sub}), the
height restriction $|\x|\leq B$ in $\E(B)$ lifts to
\begin{equation}\lab{h1}
\max\Big\{\base 1 2 2 0 1 2 3  \tau_2, ~|\tau_\ell|, ~
\base 2 3 4 3 4 5 6, ~ \base 2 2 3 1 2 3 4 |\tau_1|\Big\}\leq B,
\end{equation}
for $(\bxi,\bt)\in \TT_2$.

It will be convenient to introduce the set 
\begin{equation}
  \label{eq:FF}
\FF  :=\Big\{\bxi\in \N^7 :
\begin{array}{l}
|\mu(\xi_2\xi_3\xi_4\xi_5)|=1,\\
\hcf(\xi_1,\xi_2\xi_\ell\xi_4\xi_5)=\hcf (\xi_\ell,
\xi_2\xi_3)=1
\end{array}
\Big\}.
\end{equation}
It is easy to check that $(\bxi,\bt)\in \TT_2$ if and only if
$(\bxi,\bt)\in \FF\times \Z_*\times\N\times\Z_*$ and satisfies
\eqref{ut}, \eqref{co-2}, and \eqref{co-1'}.
With this new notation in mind, it therefore follows from Lemma \ref{bij} that
$$
\#\E(B)=
\#\Big\{
(\bxi,\bt) \in \FF 
\times \Z_*\times \N \times \Z_*:~
\text{\eqref{ut}, \eqref{co-2}, \eqref{co-1'}, \eqref{h1} hold}
\Big\}.
$$
Our final task in this section will be to analyse the inequalities
\eqref{h1} in the light of \eqref{ut}, to see
precisely what conditions they force upon the variables $\bxi,\bt$.
It is at this point of the argument that the real-valued functions
introduced in \S \ref{real} enter the picture.

Let $\bxi \in \N^7$. Then it will be convenient to define
\begin{equation}\lab{notation-1}
\al:=B^{-1/2}  \xi_1
\xi_2^{3/2}\xi_3^{2}\xi_\ell^{3/2} \xi_4^2
\xi_5^{5/2} \xi_6^3,
\end{equation}
and
\begin{equation}\lab{notation-2}
X_1:=B^{1/3}\xi_1^{-2/3}\xi_3^{-1/3}\xi_\ell
\xi_4^{2/3} \xi_5^{1/3}, \quad X_2:=
B^{1/2}\xi_2^{-1/2}\xi_\ell^{3/2} \xi_4 \xi_5^{1/2}.
\end{equation}
Now let $(\bxi,\bt)\in \TT_2$ such that \eqref{h1} holds. Then it
immediately follows from \eqref{h1} that
\begin{equation}\lab{h-al}
0\leq \al \leq 1.
\end{equation}
Recall the definition \eqref{g1} of the function $g_1:[0,1]
\rightarrow \R$.  Then it is not hard to combine \eqref{ut} and
\eqref{h1} to conclude that for fixed $\bxi$, the variable $\tau_1$
is constrained to lie in the interval
\begin{equation}\lab{h-t1}
X_1 g_1(\alpha)\leq \tau_1 \leq X_1.
\end{equation}
Next we recall the definitions \eqref{g21} and \eqref{g22} of the functions
$g_{21},g_{22}$ that are both defined on $(-\infty,1]\times [0,1]$.
For convenience, we suppose $B\in \R\smallsetminus \overline\Q$ is
such that for all $\tau_1<-X_1$ we have
$$
X_2g_{21}(\tau_1/X_1,\alpha)=X_2\sqrt{-1-(\tau_1/X_1)^3}\notin \Z.
$$
Then one deduces from \eqref{ut} and \eqref{h1} that for fixed $\bxi,
\tau_1$, the variable $\tau_2$ lies in the interval
\begin{equation}\lab{h-t2}
X_2 g_{21}(\tau_1/X_1, \alpha)< \tau_2 \leq X_2 g_{22}(\tau_1/X_1,\alpha).
\end{equation}
Putting all of this together,
and noting that we automatically have $|\tau_\ell|\leq B$
when \eqref{ut} and \eqref{h-al}--\eqref{h-t2} are satisfied, 
we have therefore established the
following result.

\begin{lem}\label{lem:count}
Let $B\in \R\smallsetminus \overline\Q$.  Then we have
\[
\#\E(B)=
\#\Big\{
(\bxi,\bt) \in \FF\times \Z_*\times \N \times \Z_*:
\begin{array}{l}
\mbox{\eqref{ut}, \eqref{co-2}, \eqref{co-1'},}\\
\mbox{\eqref{h-al}--\eqref{h-t2} hold}.
\end{array}
\Big\}.
\]
\end{lem}

\section{Proof of Theorem \ref{main'}: summations}\label{proof3}

The goal of this section is to produce a preliminary estimate
for $\#\E(B)$.  The first step, in \S \ref{2-ell}, will involve
summing over the possible values of $\tau_2,\tau_\ell$, for fixed
$\bxi,\tau_1$.  In doing so we shall take care of the summation over
$\tau_\ell$ by viewing the equation \eqref{ut} as a congruence
$$
\tau_2^2\xi_2+\tau_1^3\xi_1^2\xi_3 \eqm{0}{\xi_\ell^3\xi_4^2\xi_5}.
$$
Of course we must first carry out a  M\"obius inversion in order to
account for the coprimality condition $\hcf(\tau_\ell,\xi_4\xi_5\xi_6)=1$,
but this presents a purely technical difficulty.
The most challenging aspect involves giving a satisfactory
treatment of the error term that arises from approximating the number
of solutions $\tau_2$ to this congruence by an appropriate real-valued
function.  To facilitate digestion of
the text, the handling of the error term will be carried out separately in \S
\ref{error-10}, and this can be omitted at a first reading.  It should
be remarked that it is mainly for this part of the proof of Theorem
\ref{main'} that the contents of \S \ref{e-sum}
and \S \ref{equi} are necessary.

Next in \S \ref{1} we shall sum over the possible
values of $\tau_1$, for fixed $\bxi$.  Here again there is significant
work to be done in handling the overall contribution from the error
term. In fact, rather than showing that the overall error term makes a
negligible contribution to $\#\E(B)$, we find here that the
contribution from this error term produces a highly non-trivial
secondary contribution.  In view of the technical difficulties
involved we have therefore found it convenient to cordon off the
said treatment into \S \ref{error-11}.  This part can also be omitted at a
first reading. Finally we shall sum over the remaining
variables $\bxi$ in \S \ref{1-6}.

\subsection{Summation over  $\tau_2$ and $\tau_\ell$}\lab{2-ell}

Recall the notations \eqref{eq:FF}--\eqref{notation-2}
introduced above, and let $(\bxi,\tau_1)\in \FF \times \Z_*$ be such
that (\ref{co-1'}), (\ref{h-al})  and (\ref{h-t1}) hold.
Then the aim of this section is to estimate the number $N=N(\bxi,\tau_1)$
of $(\tau_2,\tau_\ell)\in \N\times \Z_*$ such that (\ref{co-2}),
(\ref{h-t2}) and the torsor equation (\ref{ut}) are fulfilled.
Now a M\"obius inversion yields
$$
N=\sum_{k_\ell \mid \xi_4\xi_5\xi_6}
\mu(k_\ell) N_{k_\ell},
$$
where $N_{k_\ell}$ has the same definition as $N$
except that the coprimality relation $\hcf(\tau_\ell,\xi_4\xi_5\xi_6)=1$
is removed,  and the equation (\ref{ut}) is replaced by
$$
{k_\ell}\tau_\ell\xi_\ell^3\xi_4^2\xi_5+\tau_2^2\xi_2+\tau_1^3\xi_1^2\xi_3 = 0.
$$
Let us write $I_2$ for the interval (\ref{h-t2}), and define
\begin{equation}\lab{q}
q:=k_\ell \xi_\ell^3 \xi_4^2 \xi_5.
\end{equation}
It now follows that
$$
N_{k_\ell} = \#\big\{\tau_2\in I_2\cap\N :
\hcf(\tau_2,\xi_1\xi_3)=1,
-\tau_2^2 \xi_2 \equiv \tau_1^3\xi_1^2\xi_3 \mod{q}\big\},
$$
which we henceforth assume to be non-zero.
Note that
$$
\hcf(\tau_1^3\xi_1^2\xi_3,q)=
\hcf(\tau_1^3\xi_1^2\xi_3, \xi_\ell^3\xi_4^3 \xi_5^2\xi_6, \tau_2^2\xi_2)=1,
$$
for any $k_\ell \mid \xi_4\xi_5\xi_6$. It follows from this that
\begin{equation}\lab{q123}
\hcf(\xi_1\xi_2\xi_3,k_\ell)=1.
\end{equation}
Since $k_\ell\mid q$, we must therefore take care to sum only over
values of ${k_\ell} \mid \xi_4\xi_5\xi_6$ such that
$\hcf(k_\ell, \tau_1\xi_1\xi_2\xi_3)=1$, in order to ensure a non-zero
value of $N_{k_{\ell}}$.
This latter condition is equivalent to
$\hcf(k_\ell, \xi_1\xi_2\xi_3)=1$, for any ${k_\ell} \mid
\xi_4\xi_5\xi_6$, since
$\tau_1$ is coprime to $\xi_\ell\xi_4\xi_5\xi_6$.

It is now clear that there exists a unique integer $\vr$ satisfying
$1 \leq \vr \leq q$ and $\hcf(\vr,q)=1$, such that
$$
\tau_2 \equiv \vr\tau_1\xi_1 \mod{q}, \quad -\vr^2\xi_2 \equiv \tau_1
\xi_3 \mod{q}.
$$
Thus we have
$$
N=\sum_{\colt{k_\ell \mid \xi_4\xi_5\xi_6}{
\hcf(k_\ell,\xi_1\xi_2\xi_3)=1}}
\mu(k_\ell) \sum_{\colt{1\leq \vr\leq q, ~\hcf(\vr,q)=1}{-\vr^2\xi_2\equiv
       \tau_1\xi_3 \mod{q}}}N_{k_\ell}(\vr),
$$
where
\begin{align*}
N_{k_\ell}(\vr)
&=\#\big\{\tau_2\in I_2\cap\N: \hcf(\tau_2,\xi_1\xi_3)=1, ~\tau_2
\equiv \vr\tau_1\xi_1 \mod{q}\big\}\\
&=\sum_{k_{2} \mid \xi_1\xi_3}\mu(k_2)\#\big\{\tau_2\in
(k_{2}^{-1}I_2)\cap \N :
k_2\tau_2 \equiv \vr\tau_1\xi_1 \mod{q}\big\},
\end{align*}
and in which we also know that
$\hcf(k_{2},q)\leq \hcf(\xi_1\xi_3, q) =1$.  The summand here is just the number of
positive integers restricted to lie in a certain interval, that lie in a
fixed residue class modulo $q$. When we come in the next section to sum
the error  term that arises from estimating this
quantity over all of the
remaining variables, a significant awkwardness is caused by the fact
that the interval $I_2$, as given by
(\ref{h-t2}), depends intimately upon $\tau_1$.
Write~$\ov{k_2}$
for the multiplicative inverse of $k_2$
modulo $q$.   Then it is straightforward to deduce from Lemma \ref{cong} that
\begin{align*}
N_{k_\ell}(\vr) &=
\sum_{k_2 \mid \xi_1\xi_3}
\mu(k_2)\Big(\frac{h_2(\tau_1)-h_1(\tau_1)}{q} +
r(h_1(\tau_1),h_2(\tau_1);\ov{k_2}\vr\tau_1\xi_1,q) \Big)\\
&=X_2g_2(\tau_1/X_1, \alpha)
\frac{\phi^*(\xi_1\xi_3)}{q}
+
\sum_{k_2\mid \xi_1\xi_3}
\mu(k_2) r(h_1(\tau_1),h_2(\tau_1);\ov{k_2}\vr\tau_1\xi_1,q),
\end{align*}
where we have set
\begin{equation}\lab{hi}
h_i(x):=\frac{X_2g_{2i}(x /X_1, \alpha)}{k_2},
\end{equation}
for $i=1,2$, and where $g_2$ is given by (\ref{g2'}) and
(\ref{eq:integral_gone}), and $\phi^*$ is given by~\eqref{defphi*}.
On recalling the definition (\ref{q}) of $q$, we
have therefore established the following result
for fixed $(\bxi,\tau_1)\in \FF \times \Z_*$ such that
(\ref{co-1'}), (\ref{h-al})  and~(\ref{h-t1}) hold.

\begin{lem}\lab{lem:t-2-ell}
We have
$$
N(\bxi,\tau_1)
=
\frac{X_2}{\xi_\ell^3\xi_4^2\xi_5}g_2(\tau_1/X_1,
\alpha)\Sigma(\bxi,\tau_1) + E_1(\bxi,\tau_1),
$$
where
$$
\Sigma(\bxi,\tau_1):=\phi^*(\xi_1\xi_3)
\sum_{\colt{k_\ell \mid \xi_4\xi_5\xi_6}{
\hcf(k_\ell,\xi_1\xi_2\xi_3)=1}}
\frac{\mu(k_\ell)}{k_\ell}
\sum_{\colt{1\leq \vr\leq q, ~\hcf(\vr,q)=1}{-\vr^2\xi_2\equiv
       \tau_1\xi_3 \mod{q}}}1,
$$
and
\begin{align*}
E_1(\bxi,\tau_1)&:=
\sum_{\colt{k_\ell \mid \xi_4\xi_5\xi_6}{
\hcf(k_\ell,\xi_1\xi_2\xi_3)=1}}
\sum_{k_2 \mid \xi_1\xi_3}\mu(k_2)\mu(k_\ell)\\
&\qquad
\sum_{\colt{1\leq \vr\leq q, ~\hcf(\vr,q)=1}{-\vr^2\xi_2\equiv
       \tau_1\xi_3
\mod{q}}}r(h_1(\tau_1),h_2(\tau_1);\ov{k_2}\vr\tau_1\xi_1,q).
\end{align*}
\end{lem}

For any $T \geq 1$, define the pair of sets
\begin{equation}\label{defAB}
\begin{split}
\mcal{A}(B,T)&:=\{\bxi\in \FF:
\, 0\leq \al \leq 1,\,{\xi_\ell^{3} \xi_4^2 \xi_5}\geq  T\},\\
\mcal{B}(B,T)&:=\{\bxi\in \FF:\, 0\leq\al \leq 1,
\,{\xi_\ell^{3} \xi_4^2 \xi_5}<  T\},
\end{split}
\end{equation}
where $\FF$ is given by \eqref{eq:FF} and 
$\al$ is given by \eqref{notation-1}.
As indicated above, the hardest task that we shall face is showing
that $E_1(\bxi,\tau_1)$ makes a satisfactory contribution to $\#\E(B)$,
once summed over all of the relevant values of $\bxi, \tau_1$.
That this is so is recorded in the following result, the proof of which is
postponed until the following  section.

\begin{pro}\lab{lemma10}
Let $\varepsilon>0$.  Then we have
$$
\sum_{\bxi\in \mcal{A}(B,1)}
\sum_{\colt{g_1(\al)\leq \tau_1/X_1 \leq 1}{\hcf(\tau_1,\xi_2\xi_3\xi_\ell\xi_4\xi_5\xi_6)=1}}
E_1(\bxi,\tau_1) \ll_\varepsilon
B^{43/48+\varepsilon}.
$$
\end{pro}

\subsection{Proof of Proposition \ref{lemma10}}\lab{error-10}

Recall the definition of $r(t_1,t_2;a,q)$ from the statement of Lemma
\ref{cong}, and the definitions (\ref{notation-1}), (\ref{hi}) of
$\al$ and $h_i$. Then it will suffice to estimate
\begin{align*}
Z_i(B)&:=\sum_{\bxi\in \mcal{A}(B,1)}
\sum_{\colt{k_\ell \mid \xi_4\xi_5\xi_6}{
\hcf(k_\ell,\xi_1\xi_2\xi_3)=1}}
\sum_{k_2 \mid \xi_1\xi_3}\mu(k_2)\mu(k_\ell)\\
&\qquad
\sum_{\colt{g_1(\al)\leq \tau_1/X_1 \leq 1}{
\hcf(\tau_1,\xi_2\xi_3\xi_6)=1}}
\sum_{\colt{1\leq \vr\leq q, ~\hcf(\vr,q)=1}{-\vr^2\xi_2\equiv
       \tau_1\xi_3 \mod{q}}}
\psi\Big(\frac{h_i(\tau_1)-\ov{k_2}\xi_1
\vr \tau_1}{q}\Big),
\end{align*}
for $i=1,2$.  Here, we have been able to replace the coprimality
condition $\hcf(\tau_1,\xi_2\xi_3\xi_\ell\xi_4\xi_5\xi_6)=1$ in
(\ref{co-1'}), by
$\hcf(\tau_1,\xi_2\xi_3\xi_6)=1$.
Write $Z(B)=Z_1(B)-Z_2(B)$. Then our goal is to establish the
following result, since this clearly suffices for the proof of
Proposition \ref{lemma10}.

\begin{lem}\lab{goal}
Let $\varepsilon>0$.  Then we have
$$
Z(B)=O_\varepsilon(B^{43/48+\varepsilon}).
$$
\end{lem}

Before commencing the proof of Lemma \ref{goal} proper, it will be
convenient to introduce some more notation.  Let
\begin{equation}\lab{notation-1'}
q_0:=\xi_\ell^{3} \xi_4^2 \xi_5.
\end{equation}
Then $q=k_\ell q_0$ by \eqref{q}, and
\begin{equation}\lab{notation-2'}
X_1=\frac{{q_0}^{1/3}B^{1/3}}{\xi_1^{2/3}\xi_3^{1/3}}, \quad X_2=
\frac{{q_0}^{1/2}B^{1/2}}{\xi_2^{1/2}},
\end{equation}
by \eqref{notation-2}.  It clearly follows that
\begin{equation}\lab{1/2}
\frac{X_2}{\xi_1\xi_3 X_1} =
\frac{{q_0}^{1/6}B^{1/6}}{\xi_1^{1/3}\xi_2^{1/2}\xi_3^{2/3}} \geq 1,
\end{equation}
for any $\bxi\in \N^7$ such that $0 \leq \al \leq 1$.
We also define
\begin{equation}\lab{Y12}
Y_1:=\frac{X_1}{k_1}=\frac{{q_0}^{1/3}B^{1/3}}{k_1\xi_1^{2/3}\xi_3^{1/3}},
\quad Y_2:=\frac{X_2}{k_2}=
\frac{{q_0}^{1/2}B^{1/2}}{k_2\xi_2^{1/2}},
\end{equation}
for any $k_1, k_2 \in \N$, and set
\begin{equation}\lab{fi}
f_i(x):=Y_2g_{2i}(x /Y_1, \alpha),
\end{equation}
for $i=1,2$.

We begin the proof of Lemma \ref{goal} by
applying a M\"obius inversion to remove the
coprimality
condition $\hcf(\tau_1,\xi_2\xi_3\xi_6)=1$ from the summation over
$\tau_1$. Thus we obtain
\begin{align*}
Z_i(B)&=\sum_{\bxi\in \mcal{A}(B,1)}
\sum_{\colt{k_\ell \mid \xi_4\xi_5\xi_6}{
\hcf(k_\ell,\xi_1\xi_2\xi_3)=1}} \sum_{k_2 \mid \xi_1\xi_3}
\sum_{\colt{k_1 \mid \xi_2\xi_3\xi_6}{\hcf(k_1,q)=1}}
\mu(k_1)\mu(k_2)\mu(k_\ell)\\
&\qquad
\sum_{\colt{g_1(\al)\leq \tau_1/Y_1 \leq 1}{\hcf(\tau_1,q)=1}}
\sum_{\colt{1\leq \vr\leq q}{-\vr^2\xi_2\equiv
       k_1\tau_1\xi_3 \mod{q}}}
\psi\Big(\frac{f_i(\tau_1)-k_1\ov{k_2}\xi_1
\vr \tau_1}{q}\Big),
\end{align*}
for $i=1,2$, where $f_1,f_2$ are given above.
Here we have inserted the necessary coprimality conditions
$\hcf(k_1,q)=1$ and $\hcf(\tau_1,q)=1$, and then removed
reference to the condition $\hcf(\vr,q)=1$ in the the summation over
$\vr$, since this must automatically hold.

Let $I \subset (-\infty,1]$ be an interval and
let $f$ be a real-valued function
on $I$.  Then on setting
$$
Q_I(f):=
S_I(f,q;-\xi_2,k_1\xi_3,k_1\ov{k_2}\xi_1),
$$
in the notation of \eqref{S_I},
we therefore conclude that
\begin{equation}\lab{beach}
|Z_i(B)|\leq \sum_{\bxi\in \mcal{A}(B,1)}
\sum_{\colt{k_\ell \mid \xi_4\xi_5\xi_6}{
\hcf(k_\ell,\xi_1\xi_2\xi_3)=1}} 
\sum_{k_2 \mid \xi_1\xi_3}
\sum_{\colt{k_1 \mid \xi_2\xi_3\xi_6}{\hcf(k_1,q)=1}}
|Q_{[g_1(\al)Y_1,Y_1]}(f_i)|,
\end{equation}
for $i=1,2$ .
Since
$\hcf(k_1\ov{k_2}\xi_1\xi_2\xi_3,q)=1$, the path is now clear to try
and apply Lemmas \ref{Lucca}, \ref{cor1} and \ref{cor2} to estimate
the right hand side of \eqref{beach}.

We recall the definition \eqref{defAB} of $\mcal{A}(B,T)$ and
$\mcal{B}(B,T)$. Then for $i=1,2$, we let $U_i(B,T)$ denote the overall
contribution to the right hand side of \eqref{beach} arising from
summing the $\bxi$ over the set $\mcal{A}(B,T)$, and we let
$V_i(B,T)$ denote the corresponding
contribution from the set $\mcal{B}(B,T)$.
We will then have
\begin{equation}\lab{final}
|Z(B)| \leq U_1(B,T)+U_2(B,T)+V_1(B,T)+V_2(B,T),
\end{equation}
and it will remain to choose a suitable value of $T\geq 1$ that
minimises the bound.
The hardest of the quantities to estimate is $U_i(B,T)$, for
$i=1,2$, and it is these two sums that we consider first.
It is likely that some small improvement is possible
in the exponents of these results.
However we have chosen not to pursue this here, since what we have
achieved is sufficient for our purposes.

\begin{lem}\lab{V_1}
Let $\varepsilon>0$.  Then for any $T \geq \sqrt{B}$ we have
$$
U_1(B,T) \ll_\varepsilon
B^{43/48+\varepsilon}+\frac{B^{3/2+\varepsilon}}{T^{10/9}}.
$$
\end{lem}

\begin{proof}
The thrust of the proof is concerned with estimating
$Q_{[g_1(\al)Y_1,Y_1]}(f_1).$ In fact it will be
necessary to break the interval
$[g_1(\al)Y_1,Y_1]$ into a
disjoint union of smaller intervals, on which we
can apply Lemmas \ref{cor1} and \ref{cor2} most
effectively.
This reflects the fact that for certain values of $x \in
[g_1(\al)Y_1,Y_1]$ the functions $f_1'(x)$ and
$f_1''(x)$ can be abnormally large or small.
Let 
$$\de=\de(B), \quad \eta=\eta(B), 
$$
be certain positive functions of the
parameter $B$, such that $\de >\eta$.
These functions will be selected in due course,
but we remark at the outset that $\delta < B^{-1/1000}$.
We now write
$[g_1(\al)Y_1, Y_1]=\bigcup_{i=1}^7 I_i$,  with
\begin{align*}
I_1&:=[g_1(\al)Y_1,-(2^{2/3}+\eta^{3/4})Y_1], \\
I_2&:=(-(2^{2/3}+\eta^{3/4})Y_1,-(2^{2/3}-\eta^{3/4})Y_1],\\
I_3&:=(-(2^{2/3}-\eta^{3/4})Y_1,-\mbox{$\frac{2^{2/3}+1}{2}$}Y_1], \\
I_4&:=(-\mbox{$\frac{2^{2/3}+1}{2}$}Y_1,-(1+\de)Y_1], \\
I_5&:=(-(1+\de)Y_1,-(1+\eta)Y_1], \\
I_6&:=(-(1+\eta)Y_1,-Y_1],\\
I_7&:=(-Y_1,Y_1].
\end{align*}
These intervals are clearly all non-empty, for sufficiently large
values of $B$,  and it is not hard to see that
\begin{equation}\lab{ox1}
\begin{array}{ll}
\rom{meas}(I_1)\ll Y_1/\al^{2/3}, &\quad
\rom{meas}(I_2)=2\eta^{3/4} Y_1,\\
\rom{meas}(I_3)\ll Y_1, &\quad \rom{meas}(I_4)\ll Y_1,\\
\rom{meas}(I_5)\ll \de Y_1, &\quad  \rom{meas}(I_6)=\eta Y_1,\\
\rom{meas}(I_7)=2Y_1. &
\end{array}
\end{equation}
Here we have used the definition \eqref{g1} of $g_1$, which shows in
particular that $g_1$ is negative, to deduce that
\begin{equation}\lab{tower}
g_1(\al) \ll \al^{-2/3}.
\end{equation}
Before proceeding with the task of estimating $Q_i:=Q_{I_i}(f_1)$ for $1
\leq i \leq 7$, we first deduce from  \eqref{g21} and
\eqref{fi} that
\begin{equation}\lab{siena1}
f_1(x)=
\left\{
\begin{array}{ll}
Y_2\sqrt{-1-(x /Y_1)^3}, &
x \in I_1\cup\cdots \cup I_6,\\
0, & x \in I_7.
\end{array}
\right.
\end{equation}
In particular it follows that on the interior of
each interval $I_i$, the real-valued
function $f_1$ is infinitely differentiable, with $f_1'$ being
monotonic and of constant sign.

Moreover, we see that
$$
\Big|\int_{I_i}f_1'(t)\d t\Big|+1 \ll \la_0^{(i)}
$$
for $1\leq i \leq 7$, with
\begin{equation}\lab{ox2}
\begin{array}{lll}
\la_0^{(1)}= Y_2/\al, &\quad
\la_0^{(2)}= \eta^{3/4} Y_2 + 1, &\quad
\la_0^{(3)}=\la_0^{(4)}= Y_2,\\
\la_0^{(5)}= \de^{1/2} Y_2+1, &\quad
\la_0^{(6)}= \eta^{1/2} Y_2 +1, &\quad
\la_0^{(7)}=1.
\end{array}
\end{equation}
It follows from these remarks, together with \eqref{def-sig},
that $f_1 \in C^1(I_i;A\la_0^{(i)})$ for
$1 \leq i\leq 7$, for some absolute constant $A>0.$
%t%

We begin by using Lemma \ref{cor1} to estimate $Q_2, Q_6$ and $Q_7$.
Thus it follows from the trivial bounds $Y_1\leq X_1$ and $Y_2\leq
X_2$, that
\begin{equation}\lab{H267}
\begin{split}
\sum_{i=2,6,7}|Q_i|
&\ll_\varepsilon
q^{\varepsilon}X_1^\varepsilon\Big(
q^{1/2}+
\frac{\eta^{3/4}X_1^{1/2}X_2^{1/2}}{q^{1/4}}+
\frac{X_1+X_2^{1/2}}{q^{1/4}}\Big).
\end{split}
\end{equation}
Turning to the size of $Q_1, Q_3,Q_4$ and $Q_5$,
we must examine in more detail the behaviour of the function $f_1$
on $I_1\cup I_3\cup I_4\cup I_5$.   Now it is clear that
$$
f_1'(x)= \frac{-3Y_2(x/Y_1)^2}{2Y_1\sqrt{-1-(x/Y_1)^3}}, \quad
f_1''(x)=
\frac{3Y_2(x/Y_1)\Big(4+(x/Y_1)^3\Big)}{4Y_1^2(-1-(x/Y_1)^3)^{3/2}}.
$$
In particular $f_1'$ has constant sign on
$I_1\cup I_3\cup I_4\cup I_5$, and is monotonic
increasing (resp. decreasing)
on $I_1$ (resp. on $I_3\cup I_4\cup I_5$). We
therefore deduce from \eqref{tower} that
$$
f_1'(x)\ll \left\{
\begin{array}{ll}
\al^{-4/3}Y_2/Y_1, & x \in I_{1},\\
Y_2/Y_1, & x \in I_{3},\\
\de^{-1/2}Y_2/Y_1, & x \in I_{4},\\
\eta^{-1/2}Y_2/Y_1, & x \in I_{5}.
\end{array}
\right.
$$
Moreover, it is clear from \eqref{1/2} and \eqref{Y12} that
$$
Y_2/Y_1 \geq X_2/(\xi_1\xi_3 X_1) \geq 1.
$$
On $(-\infty,-Y_1]$ it is clear that the function
$|f_1''(x)|$ is minimised at $x=-2^{2/3}Y_1$, where it takes the value
$0$.  Hence on $I_1\cup I_3$ we have
$$
\eta^{3/4} Y_2/Y_1^2 \ll |f_1''(x)| \ll
Y_2/Y_1^2,
$$
on $I_4$ we have
$$
Y_2/Y_1^2 \ll |f_1''(x)| \ll  \de^{-3/2} Y_2/Y_1^2,
$$
and finally on $I_5$ we have
$$
\de^{-3/2} Y_2/Y_1^2 \ll |f_1''(x)| \ll  \eta^{-3/2} Y_2/Y_1^2.
$$
On recalling \eqref{ox2}, we have therefore shown that there exists an
absolute constant $A>0$ such that $f_1$ belongs to the sets
\begin{align*}%t%
&C^2\Big(I_1;A\la_0^{(1)},\frac{AY_2}{\al^{4/3}Y_1},
\frac{c_1\eta^{3/4} Y_2}{Y_1^2},
\frac{1}{\eta^{3/4}}\Big), \,
C^2\Big(I_3;A\la_0^{(3)},\frac{AY_2}{Y_1}, \frac{c_3\eta^{3/4} Y_2}{Y_1^2},
\frac{1}{\eta^{3/4}}\Big),\\
&C^2\Big(I_4;A\la_0^{(4)},\frac{AY_2}{\de^{1/2}Y_1}, \frac{c_4 Y_2}{Y_1^2},
\frac{1}{\de^{3/2}}\Big),\,
C^2\Big(I_5;A\la_0^{(5)},\frac{AY_2}{\eta^{1/2}Y_1},
\frac{c_5 Y_2}{\de^{3/2}Y_1^2},
\frac{\de^{3/2}}{\eta^{3/2}}\Big),
\end{align*}
for appropriate absolute constants $c_1,c_3,c_4,c_5>0$.
Let us write
$$
f_1 \in C^2(I_i;A\la_0^{(i)},\la_1^{(i)},\la_2^{(i)},j^{(i)}) %t%
$$
for $i=1,3,4,5$, and recall the estimates \eqref{ox1} for $\rom{meas}(I_i)$.
We are now ready to complete our estimates for $Q_1,Q_3,Q_4,Q_5$
via an application of Lemma~\ref{cor2}.
Thus we obtain
\begin{align*}
Q_1
&\ll_\varepsilon
q^{\varepsilon}Y_1^\varepsilon
\Big(
q^{1/2}+ \frac{Y_1}{\al^{2/3}q^{1/3}}+
\frac{Y_1^{1/3}Y_2^{1/3}}{\al^{10/9}\eta^{1/4}}
+\frac{\eta^{1/2}Y_1^{1/3}Y_2^{4/3}}{\al^{13/9}q^{3/2}}+
\frac{Y_2}{\al^{4/3} q}\Big),\\
Q_3
&\ll_\varepsilon
q^{\varepsilon}Y_1^\varepsilon\Big(
q^{1/2}+ \frac{Y_1}{q^{1/3}}+
\frac{Y_1^{1/3}Y_2^{1/3}}{\eta^{1/4}}
+\frac{\eta^{1/2}Y_1^{1/3}Y_2^{4/3}}{q^{3/2}}+
\frac{Y_2}{q}\Big)\\
Q_4
&\ll_\varepsilon
q^{\varepsilon}Y_1^\varepsilon\Big(
q^{1/2}+ \frac{Y_1}{q^{1/3}}+
\frac{Y_1^{1/3}Y_2^{1/3}}{\de^{1/3}}
+\frac{\de^{1/6}Y_1^{1/3}Y_2^{4/3}}{q^{3/2}}+
\frac{Y_2}{\de^{3/2} q}\Big),
\end{align*}
and
\begin{align*}
Q_5
&\ll_\varepsilon
q^{\varepsilon}Y_1^\varepsilon\Big(
q^{1/2}+ \frac{\de Y_1}{q^{1/3}}+
\frac{\de^{5/6}Y_1^{1/3}Y_2^{1/3}}{\eta^{1/3}}
+\frac{\de^{5/6}\eta^{1/6}Y_1^{1/3}Y_2^{4/3}}{q^{3/2}}+
\frac{\de^2 Y_2}{\eta^{3/2} q}\Big).
\end{align*}
It turns out that we shall need to minimise the third terms in these
estimates for $Q_1, Q_3, Q_4$ and $Q_5$.  We therefore make the selection
$$
\delta=\eta^{2/7},
$$
so that in particular $\de>\eta$.  With this choice we may conclude that
\begin{equation}\lab{H1345}
\begin{split}
\sum_{i=1,3,4,5} |Q_i|
\ll_\varepsilon&
q^{\varepsilon}X_1^\varepsilon
\Big(
q^{1/2}+ \frac{X_1}{\al^{2/3}q^{1/3}}+
\frac{X_1^{1/3}X_2^{1/3}}{\al^{10/9}\eta^{1/4}}
+\frac{X_1^{1/3}X_2^{4/3}}{\al^{13/9}q^{3/2}}\\
&\quad+\frac{X_2}{\al^{4/3} q} + \frac{X_2}{\eta q}\Big),
\end{split}
\end{equation}
since $Y_1 \leq X_1, Y_2 \leq X_2$ and $\al, \eta\leq 1$.

Write $J$ for the interval $[g_1(\al)Y_1,Y_1]$.  Then
\eqref{H267} and \eqref{H1345}  provide us with an overall
estimate for $Q_{J}(f_1)$.  This estimate is made somewhat complicated
by the large number of terms that it contains, and so we proceed to
consider the overall contribution to the right hand side of
\eqref{beach} from some of the individual terms.
Let us begin by handling the term
$q^{1/2}$ that appears in both \eqref{H267} and \eqref{H1345}.
 From \eqref{q} and the inequality $k_\ell \leq \xi_4\xi_5\xi_6$ it
follows that $q \leq \xi_\ell^3\xi_4^3\xi_5^2\xi_6$.  Thus we see that
\begin{equation}\lab{term1}
\begin{split}
\sum_{\xi_1^2\xi_2^3\xi_3^4\xi_\ell^3 \xi_4^4 \xi_5^5 \xi_6^6 \leq B} q^{1/2}
&\leq
\sum_{\xi_1^2\xi_2^3\xi_3^4\xi_\ell^3 \xi_4^4 \xi_5^5 \xi_6^6 \leq  B}
\xi_\ell^{3/2}\xi_4^{3/2}\xi_5\xi_6^{1/2}\\
&\leq
\sum_{\xi_2^3\xi_3^4\xi_\ell^3 \xi_4^4 \xi_5^5 \xi_6^6 \leq B}
\frac{B^{1/2}}{\xi_2^{3/2}\xi_3^2\xi_4^{1/2} \xi_5^{3/2}
\xi_6^{5/2}} \ll B^{5/6},
\end{split}
\end{equation}
where the summations are over all $\bxi \in \N^7$ in the proscribed range.
Next we recall the definitions
\eqref{notation-1}, \eqref{notation-1'},
\eqref{notation-2'} of
$\al,{q_0}$ and $X_1,X_2$, together with the bounds $0\leq \al \leq 1$ and
$q\geq {q_0}>T$ that hold for any $\bxi \in \mcal{A}(B,T)$.
Then it follows that
$$
\al=B^{-1/2}\xi_1 \xi_2^{3/2}\xi_3^{2}\xi_4 \xi_5^{2} \xi_6^3{q_0}^{1/2},
$$
whence
$$
\max\Big\{\frac{X_1}{q^{1/4}}, \frac{X_1}{\al^{2/3}q^{1/3}}\Big\} \leq
\frac{X_1}{\al^{2/3}{q_0}^{1/4}} \leq \frac{B^{2/3}}{\xi_1^{4/3}
    \xi_2\xi_3^{5/3}\xi_\ell^{3/4}\xi_4^{7/6} \xi_5^{19/12}\xi_6^2}.
$$
Thus we have
\begin{equation}\lab{term2}
\sum_{\xi_1^2\xi_2^3\xi_3^4\xi_\ell^3 \xi_4^4 \xi_5^5 \xi_6^6 \leq B}
\hspace{-0.3cm}
\max\Big\{\frac{X_1}{q^{1/4}}, \frac{X_1}{\al^{2/3}q^{1/3}}\Big\}
\ll
\sum_{\xi_2 \xi_\ell \leq  B^{1/3}}
\frac{B^{2/3}}{\xi_2\xi_\ell^{3/4}}\ll B^{3/4}.
\end{equation}
Similarly one notes that
$$
\max\Big\{
\frac{X_2^{1/2}}{q^{1/4}},
\frac{X_2}{\al^{4/3}q},
\frac{X_2}{\eta q}
\Big\} \leq
\frac{B^{1/4}}{\xi_2^{1/4}}
+
\frac{B^{1/2}}{\eta \xi_2^{1/2}{q_0}^{1/2}}
+
\frac{B^{7/6}}{\xi_1^{4/3}
\xi_2^{5/2}\xi_3^{8/3}\xi_\ell \xi_4^{2} \xi_5^{3}\xi_6^4{q_0}^{2/3}},
$$
whence
\begin{equation}\lab{term3}
\begin{split}
\sum_{\xi_1^2\xi_2^3\xi_3^4\xi_\ell^3 \xi_4^4 \xi_5^5 \xi_6^6 \leq B}
\hspace{-0.3cm}
\max\Big\{
\frac{X_2^{1/2}}{q^{1/4}},
\frac{X_2}{\al^{4/3}q},
\frac{X_2}{\eta q}
\Big\}
&\ll \frac{B^{3/4}}{\eta}+B^{5/6}\log B.
\end{split}
\end{equation}
Here we have used the fact that ${q_0}>T\geq \sqrt{B}$.
Now let us treat the term
$$
\frac{X_1^{1/3}X_2^{4/3}}{\al^{\theta}q^{3/2}}
\leq
\frac{B^{7/9+\theta/2}}{\xi_\ell(\xi_1\xi_2\xi_3\xi_4
\xi_5\xi_6)^{\theta}{q_0}^{7/18+\theta/2}},
$$
for any $\theta >1$.
But for any such $\theta$ we see that
\begin{equation}\lab{term4}
\sum_{\xi_1^2\xi_2^3\xi_3^4\xi_\ell^3 \xi_4^4 \xi_5^5 \xi_6^6 \leq B}
\frac{X_1^{1/3}X_2^{4/3}}{\al^\theta q^{3/2}}
\ll \frac{B^{7/9+\theta/2}\log B}{T^{7/18+\theta/2}},
\end{equation}
since ${q_0}>T$.

It turns out that for the choice of $\eta$ and $T$ that we shall make,
all of the terms \eqref{term1}--\eqref{term4} make a negligible
contribution in our final estimate for $U_1(B,T)$.  There are in
effect two dominant contributions, the first of which is the term
$$
\frac{\eta^{3/4}X_1^{1/2}X_2^{1/2}}{q^{1/4}}\leq
\frac{\eta^{3/4}{q_0}^{1/6}B^{5/12}}{\xi_1^{1/3}\xi_2^{1/4}\xi_3^{1/6}}=
\frac{\eta^{3/4}B^{5/12}\xi_\ell^{1/2}\xi_4^{1/3}\xi_5^{1/6}}{\xi_1^{1/3}\xi_2^{1/4}\xi_3^{1/6}},
$$
in \eqref{H267}, which therefore produces a contribution
\begin{equation}\lab{term5}
\sum_{\xi_1^2\xi_2^3\xi_3^4\xi_\ell^3 \xi_4^4 \xi_5^5 \xi_6^6 \leq B}
\hspace{-0.3cm}
\frac{\eta^{3/4}X_1^{1/2}X_2^{1/2}}{q^{1/4}}
\ll \eta^{3/4}B^{11/12}.
\end{equation}
The second major contribution comes from the term
$$
\frac{X_1^{1/3}X_2^{1/3}}{\al^{10/9}\eta^{1/4}}\leq
\frac{B^{5/6}}{\eta^{1/4}\xi_\ell^{5/6}(\xi_1\xi_2\xi_3\xi_4\xi_5\xi_6)^{10/9}}, 
$$
in \eqref{H1345}, which plainly yields the contribution
\begin{equation}\lab{term7}
\sum_{\xi_1^2\xi_2^3\xi_3^4\xi_\ell^3 \xi_4^4 \xi_5^5 \xi_6^6 \leq B}
\frac{X_1^{1/3}X_2^{1/3}}{\al^{10/9}\eta^{1/4}}\ll \frac{B^{8/9}}{\eta^{1/4}}.
\end{equation}

Bringing together \eqref{term1}--\eqref{term7} in the right-hand
side of \eqref{beach}, we may therefore conclude that
$$
U_1(B,T)
\ll_\varepsilon
\eta^{3/4}B^{11/12+\varepsilon}+\frac{B^{8/9+\varepsilon}}{\eta^{1/4}}+
\frac{B^{3/4+\varepsilon}}{\eta}+
\frac{B^{3/2+\varepsilon}}{T^{10/9}},
$$
for any $\varepsilon>0$.  The first two terms represent the dominant
contribution, and so it remains to choose a value of $\eta$ that
balances them.  The proof of Lemma \ref{V_1} is thus completed by taking
$\eta=B^{-1/36}$.
\end{proof}

\begin{lem}\lab{V_2}
Let $\varepsilon>0$.  Then for any $T \geq \sqrt{B}$ we have
$$
U_2(B,T) \ll_\varepsilon
B^{59/66+\varepsilon}+\frac{B^{5/3+\varepsilon}}{T^{23/18}}.
$$
\end{lem}

\begin{proof}
The proof of this result is very similar to that of the previous
lemma, and so we shall be brief.
Now for any $\al >0$ it is easy to see that
$$
g_1(\al)\leq -(1/\al^2-1)^{1/3}.
$$
As above it will be necessary to break the interval $[g_1(\al)Y_1,Y_1]$ into a
disjoint union of smaller intervals, in order to estimate
$Q_{[g_1(\al)Y_1,Y_1]}(f_2)$ most effectively.
Let $\eta$ be a small positive real number. This will be selected in due
course, but it may be assumed that $\eta < B^{-1/1000}$.
Thus we write $[g_1(\al)Y_1,Y_1]=\bigcup_{i=1}^5 I_i$,  with
\begin{align*}
I_1&:=[g_1(\al)Y_1,-(1/\al^2-1)^{1/3}Y_1],\\
I_2&:=(-(1/\al^2-1)^{1/3}Y_1,-\eta Y_1],\\
I_3&:=(-\eta Y_1,\eta Y_1], \\
I_4&:=(\eta Y_1,(1-\eta^{8/3})Y_1],\\
I_5&:=((1-\eta^{8/3}) Y_1,Y_1].
\end{align*}
Observe that
\begin{equation}\lab{ox3}
\begin{array}{lll}
\rom{meas}(I_1)\ll Y_1, &\quad
\rom{meas}(I_2)\ll Y_1/\al^{2/3}, &\quad
\rom{meas}(I_3)=2\eta Y_1,\\
\rom{meas}(I_4)\ll Y_1, &\quad
\rom{meas}(I_5)=\eta^{8/3} Y_1. &
\end{array}
\end{equation}
Before proceeding with the task of estimating $Q_i:=Q_{I_i}(f_2)$ for $1
\leq i \leq 5$, we first deduce from  \eqref{g22} and \eqref{fi} that
$$
f_2(x)=
\left\{
\begin{array}{ll}
\al^{-1}Y_2, & x \in I_1,\\
Y_2\sqrt{1-(x /Y_1)^3}, &
x \in I_2\cup I_3\cup I_4\cup I_5.
\end{array}
\right.
$$
In particular it follows that on each interval $I_i$, the real-valued
function $f_2$ is infinitely differentiable, with $f_2'$ being
monotonic and of constant sign.  Moreover,
$$
\Big|\int_{I_i}f_2'(t)\d t\Big|+1 \ll \la_0^{(i)}
$$
for $1\leq i \leq 5$, with
\begin{equation}\lab{ox4}
\begin{array}{lll}
\la_0^{(1)}=1, &\quad
\la_0^{(2)}= Y_2/\al, &\quad
\la_0^{(3)}= \eta^3 Y_2+1,\\
\la_0^{(4)}= Y_2, &\quad
\la_0^{(5)}= \eta^{4/3} Y_2 +1. &
\end{array}
\end{equation}
Hence $f_2 \in C^1(I_i;A\la_0^{(i)})$ for $1 \leq i\leq 5$, for an
absolute constant $A>0$. %t%

We begin by using Lemma~\ref{cor1} to estimate $Q_1,Q_3$ and $Q_5$.
Thus it follows from \eqref{ox3} and \eqref{ox4} that
\begin{equation}\lab{F13}
\begin{split}
\sum_{i=1,3,5}|Q_i| &\ll_\varepsilon
q^\varepsilon Y_1^\varepsilon\Big(
q^{1/2} + \frac{Y_1}{q^{1/3}}+\frac{\eta^2 Y_1^{1/2}Y_2^{1/2}}{q^{1/4}}
+\frac{Y_2^{1/2}}{q^{1/4}}
\Big).
\end{split}
\end{equation}
Turning to the size of $Q_{2}$ and $Q_{4}$, we
must examine the behaviour of $f_2$
on $I_2\cup I_4$.   Now it is clear that
$$
f_2'(x)= \frac{-3Y_2(x/Y_1)^2}{2Y_1\sqrt{1-(x/Y_1)^3}}, \quad
f_2''(x)=
\frac{-3Y_2(x/Y_1)\Big((x/Y_1)^3-4\Big)}{4Y_1^2(1-(x/Y_1)^3)^{3/2}}.
$$
In particular $f_2'$ has constant sign on
$I_2$, and $|f_2'|$ is monotonic decreasing (resp. increasing)
on $I_2$ (resp. on $I_4$). We therefore deduce that
$$
f_2'(x)\ll \left\{
\begin{array}{ll}
\al^{-1/3}Y_2/Y_1, & x \in I_{2},\\
\eta^{-4/3}Y_2/Y_1, & x \in I_{4}.
\end{array}
\right.
$$
Next we note that on $(-\infty,Y_1]$ the function
$|f_1''(x)|$ is minimised at $x=0$, where it takes the value
$0$.  Hence we have
$$
\eta Y_2/Y_1^2 \ll |f_2''(x)| \ll
\left\{
\begin{array}{ll}
Y_2/Y_1^2, & \mbox{if $x \in I_2$},\\
\eta^{-4}Y_2/Y_1^2, & \mbox{if $x \in I_4$}.
\end{array}\right.
$$
We have therefore shown that there exists an absolute constant $A>0$
such that %t%
$$
f_2  \in
C^2\Big(I_2;A\la_0^{(2)},\frac{AY_2}{\al^{1/3}Y_1},
\frac{c_2\eta Y_2}{Y_1^2},
\frac{1}{\eta}\Big), \quad
f_2  \in C^2\Big(I_4;A\la_0^{(4)},\frac{AY_2}{\eta^{4/3}Y_1}, \frac{c_4\eta Y_2}{Y_1^2},
\frac{1}{\eta^{5}}\Big),
$$
for appropriate absolute constants $c_2,c_4>0$.
Let us write
$$
f_2 \in C^2(I_i;A\la_0^{(i)},\la_1^{(i)},\la_2^{(i)},j^{(i)}) %t%
$$
for $i=2,4$. On recalling \eqref{ox3}, we may
therefore deduce from Lemma \ref{cor2} that
\begin{equation}\lab{F2}
Q_{2}
\ll_\varepsilon
q^{\varepsilon}Y_1^\varepsilon\Big(
q^{1/2}+ \frac{Y_1}{\al^{2/3}q^{1/3}}+
\frac{Y_1^{1/3}Y_2^{1/3}}{\al^{4/9}\eta^{1/3}}
+\frac{\eta^{2/3}Y_1^{1/3}Y_2^{4/3}}{\al^{16/9}q^{3/2}}+
\frac{Y_2}{\al^{4/3}q}\Big),
\end{equation}
and
\begin{equation}\lab{F4}
Q_{4}
\ll_\varepsilon
q^{\varepsilon}Y_1^\varepsilon\Big(
q^{1/2}+ \frac{Y_1}{q^{1/3}}+
\frac{Y_1^{1/3}Y_2^{1/3}}{\eta^{4/9}}
+\frac{\eta^{10/9}Y_1^{1/3}Y_2^{4/3}}{q^{3/2}}+
\frac{Y_2}{\eta^{4}q}\Big).
\end{equation}

Write $J$ for the interval $[g_1(\al)Y_1,Y_1]$.
Then the estimates \eqref{F13}--\eqref{F4} together provide us with an overall
estimate for $Q_{J}(f_2)$. Once taken into the right hand side of
\eqref{beach}, the summation over $\bxi$ being only over $\bxi \in
\mcal{A}(B,T)$, we may draw upon the proof of Lemma \ref{V_1} to
handle the contribution from all of the individual terms.
Thus it follows from \eqref{term1}--\eqref{term3} that
\begin{equation}\lab{term1'}
\sum_{\bxi \in \mcal{A}(B,T)}\Big(q^{1/2}+\frac{X_1}{\al^{2/3}q^{1/3}}
+\frac{X_2^{1/2}}{q^{1/4}}
+\frac{X_2}{\al^{4/3}q}
+\frac{X_2}{\eta^{4}q}
\Big)
\ll B^{5/6}\log B + \frac{B^{3/4}}{\eta^4},
\end{equation}
since $T \geq \sqrt{B}$.
Similarly, \eqref{term4} gives
\begin{equation}\lab{term2'}
\sum_{\bxi\in \mcal{A}(B,T)}
\frac{X_1^{1/3}X_2^{4/3}}{\al^{16/9}q^{3/2}}
\ll \frac{B^{5/3}\log B}{T^{23/18}}.
\end{equation}
These will all make a negligible contribution in our final estimate
for the quantity $U_2(B,T)$.  There are now three dominant
contributions, the first of which is the term
$\eta^{2}X_1^{1/2}X_2^{1/2}/q^{1/4}$ in
\eqref{F13}. This produces an overall contribution
\begin{equation}\lab{term3'}
\sum_{\bxi\in \mcal{A}(B,T)}
\frac{\eta^{2}X_1^{1/2}X_2^{1/2}}{q^{1/4}}
\ll \eta^{2}B^{11/12},
\end{equation}
by \eqref{term5}.
Next we note that
$$
\max\Big\{
\frac{X_1^{1/3}X_2^{1/3}}{\al^{4/9}\eta^{1/3}},
\frac{X_1^{1/3}X_2^{1/3}}{\eta^{4/9}}
\Big\} \leq \frac{X_1^{1/3}X_2^{1/3}}{\al^{4/9}\eta^{4/9}} \leq
\frac{{q_0}^{1/18}B^{1/2}}{\eta^{4/9}\xi_1^{2/3}\xi_2^{5/6}\xi_3\xi_4^{4/9}\xi_5^{8/9}\xi_6^{4/3}}.
$$
Hence we may argue as for \eqref{term7} to deduce
that the second major contribution is
\begin{equation}\lab{term4'}
\sum_{\bxi\in \mcal{A}(B,T)}
\Big(
\frac{X_1^{1/3}X_2^{1/3}}{\al^{4/9}\eta^{1/3}}
+\frac{X_1^{1/3}X_2^{1/3}}{\eta^{4/9}}
\Big) \ll
\frac{B^{8/9}}{\eta^{4/9}}.
\end{equation}

We may now bring together \eqref{term1'}--\eqref{term4'} in the right-hand
side of \eqref{beach}, in order to conclude that
$$
U_2(B,T)
\ll_\varepsilon
\eta^{2}B^{11/12+\varepsilon}+\frac{B^{8/9+\varepsilon}}{\eta^{4/9}}+
\frac{B^{3/4+\varepsilon}}{\eta^4}+
\frac{B^{5/3+\varepsilon}}{T^{23/18}},
$$
for any $\varepsilon>0$.  As in the proof of
Lemma \ref{V_1}, the first two terms represent
the dominant contribution and we therefore select $\eta=B^{-1/88}$. This completes
the proof of Lemma \ref{V_2}.
\end{proof}

Our final task is to estimate the sizes of $V_1(B,T)$ and $V_2(B,T)$.
This is entirely straightforward, since it will suffice just to employ
a trivial upper bound for the sums $Q_{[g_1(\al)Y_1,Y_1]}(f_1)$ and
$Q_{[g_1(\al)Y_1,Y_1]}(f_2)$.  Thus it follows from \eqref{tower} and
Lemma \ref{Lucca} that
$$
Q_{[g_1(\al)Y_1,Y_1]}(f_i)\ll_\varepsilon \frac{q^\varepsilon Y_1}{\al^{2/3}}
\ll_\varepsilon
\frac{B^{2/3+\varepsilon}}{\xi_1^{4/3}\xi_2\xi_3^{5/3}\xi_4^{2/3}\xi_5^{4/3}\xi_6^{2}},
$$
on substituting the definitions of $\al, Y_1$ into this estimate.
But by definition of the set $\mcal{B}(B,T)$, we
must have ${q_0}\leq T$, whence
$$
\xi_\ell \leq\frac{T^{1/3}}{\xi_4^{2/3}\xi_5^{1/3}}.
$$
On summing over all values of $\bxi \in\mcal{B}(B,T)$ we therefore
obtain the following result.

\begin{lem}\lab{W_12}
Let $\varepsilon>0$ and let $i=1$ or $2$.  Then for any $T \geq 1$ we have
$$
V_i(B,T) \ll_\varepsilon T^{1/3}B^{2/3+\varepsilon}.
$$
\end{lem}

We are now in a position to draw together Lemmas \ref{V_1}, \ref{V_2}
and \ref{W_12} in \eqref{final}.  Thus we conclude that
$$
Z(B) \ll_\varepsilon
B^{43/48+\varepsilon}+\frac{B^{3/2+\varepsilon}}{T^{10/9}}+
\frac{B^{5/3+\varepsilon}}{T^{23/18}}+
T^{1/3}B^{2/3+\varepsilon},
$$
for any $\varepsilon>0$ and any $T \geq \sqrt{B}$.
We therefore complete the proof of Lemma \ref{goal}, and so the proof of
Proposition \ref{lemma10},  by taking
$T =B^{18/29}$.

\subsection{Summation over $\tau_1$}\lab{1}

For fixed $\bxi \in \FF$ such that 
\eqref{h-al} holds, we proceed to sum the main term in Lemma
\ref{lem:t-2-ell} over all $\tau_1$ satisfying \eqref{co-1'} and
\eqref{h-t1}. Thus our task is to estimate 
\[N' = N'(\bxi) :=
\frac{X_2}{\xi_\ell^3\xi_4^2\xi_5}
\sum_{\colt{g_1(\al)\leq \tau_1/X_1 \leq
1}{\hcf(\tau_1,\xi_2\xi_3\xi_\ell\xi_4\xi_5\xi_6)=1}}
g_2(\tau_1/X_1, \alpha)\Sigma(\bxi,\tau_1),
\]
where $g_2$ is given by \eqref{g2'} and $\Sigma(\bxi,\tau_1)$ is
as in the statement of Lemma \ref{lem:t-2-ell}.

Let $t_1,t_2 \in \R$ such that $t_2 \geq  t_1$.
Then we begin by deriving an asymptotic formula for
\begin{equation}\lab{bike}
\NN(t_1,t_2) :=
\phi^*(\xi_1\xi_3)
\sum_{\colt{k_\ell \mid \xi_4\xi_5\xi_6}{
\hcf(k_\ell,\xi_1\xi_2\xi_3)=1}} 
\frac{\mu(k_\ell)}{k_\ell}
\sum_{\colt{1\leq \vr\leq q}{\hcf(\vr,q)=1}} N'_{k_\ell}(\varrho;t_1,t_2),
\end{equation}
where $q$ is given by \eqref{q} and
\[
N'_{k_\ell}(\varrho;t_1,t_2) = \#\Big\{\tau_1 \in [t_1,t_2]:
\begin{aligned}
{} & \hcf (\tau_1,\xi_2\xi_3\xi_\ell\xi_4\xi_5\xi_6)=1,
\\ &
-\vr^2\xi_2\equiv \tau_1\xi_3 \mod{q}
\end{aligned}
\Big\}.
\]
In view of the fact that $\hcf(\varrho^2\xi_2, q) = 1$, we may plainly
replace the condition
$\hcf(\tau_1,\xi_2\xi_3\xi_\ell\xi_4\xi_5\xi_6)=1$
by $\hcf(\tau_1,\xi_2\xi_3\xi_6) = 1$
in the definition of $N'_{k_\ell}(\varrho;t_1,t_2)$.  On performing
a further M\"obius inversion, we obtain
\begin{align*}
\NN(t_1,t_2) &=
\phi^*(\xi_1\xi_3)
\hspace{-0.3cm}
\sum_{\colt{k_\ell \mid \xi_4\xi_5\xi_6}{
\hcf(k_\ell,\xi_1\xi_2\xi_3)=1}} 
\hspace{-0.2cm}
\frac{\mu(k_\ell)}{k_\ell}
\hspace{-0.2cm}
   \sum_{\colt{k_1\mid \xi_2\xi_3\xi_6}{\hcf(k_1,q)=1}}
\hspace{-0.1cm}
\mu(k_1)
\hspace{-0.1cm}
\sum_{\colt{1\leq \vr\leq q}{\hcf(\vr,q)=1}}
\hspace{-0.1cm}
N'_{k_1,{k_\ell}}(\varrho; t_1,t_2),
\end{align*}
where
\[
N'_{k_1,{k_\ell}}(\varrho;t_1,t_2) =
\#\big\{\tau_1 \in [t_1/k_1, t_2/k_1]:
\congr{-\varrho^2\xi_2}{k_1\tau_1\xi_3}{q} \big\}.
\]
We have therefore reduced the problem to once again estimating the
number of integers which are restricted to lie in a
certain interval, and that lie in a fixed residue class modulo $q$.

Let $a = a(\bxi,k_1,{k_\ell})$ be the unique positive integer such that $a
\le q$, with $\hcf(a,q)=1$ and
\[
\congr{-\xi_2}{k_1 a \xi_3}{q}.
\]
Then $\congr{-\varrho^2\xi_2}{k_1 \tau_1\xi_3}{q}$ if and only if
$\congr{\tau_1}{a \varrho^2}{q}$. Using
Lemma~\ref{cong}, we therefore conclude that
\[
N'_{k_1,{k_\ell}}(\varrho;t_1,t_2) = \frac{t_2-t_1}{k_1q} +
r(t_1/k_1,t_2/k_1;a\varrho^2,q),
\]
if $t_1\notin \Z$.
Now let
\begin{equation}\lab{vt}
   \vartheta(\bxi) =
\phi^*(\xi_2\xi_3\xi_\ell\xi_4\xi_5\xi_6)
\frac{\phi^*(\xi_4\xi_5\xi_6)\phi^*(\xi_1\xi_3)}{\phi^*(\hcf(
\xi_6,\xi_1\xi_2\xi_3))},
\end{equation}
if $\bxi\in\FF$, and $\vt(\bxi)=0$ otherwise. Then for any $t_1,t_2 \in \R$
such that $t_2\geq t_1$ and 
$t_1\notin\Z$, we have therefore shown that
\[
\NN(t_1,t_2) =
\vt(\bxi)(t_2-t_1) + \cR(t_1,t_2)
\]
in (\ref{bike}), where
\begin{equation}\lab{bike2}
\begin{split}
\cR(t_1,t_2)
&=\phi^*(\xi_1\xi_3)\sum_{\colt{k_\ell\mid
\xi_4\xi_5\xi_6}{
\hcf(k_\ell,\xi_1\xi_2\xi_3)=1}} 
{\mu(k_\ell)\over
k_\ell}\sum_{\colt{k_1\mid
\xi_2\xi_3\xi_6}{\hcf(k_1,q)=1}}\mu(k_1)\\
&\qquad \sum_{\colt{1\leq \vr\leq q}{\hcf(\vr,q)=1}}
r(t_1/k_1,t_2/k_1;a\vr^2,q).
\end{split}
\end{equation}
Define the function $g_3:[0,1]\rightarrow \R$ on the unit interval, given by
\begin{equation}\label{g3}
\begin{split}
g_3(v) &:= \int_{g_1(v)}^1 {g_2}(u,v) \dd u.
\end{split}
\end{equation}
Then a straightforward application of partial summation yields
the following result, in  which $D_1{g_2}$ is the derivative of
${g_2}$ with respect to the first variable, and
$\bxi \in \FF$ is such that \eqref{h-al} holds.

\begin{lem}\label{lem:t-1}
We have
\[
N'(\bxi)=
\frac{\vartheta(\bxi)X_1X_2}{\xi_\ell^3\xi_4^2\xi_5} g_3(\al) + E_2(\bxi),
\]
with
\[
E_2(\bxi):=
\frac{-X_2}{\xi_\ell^3\xi_4^2\xi_5}
\int_{-\al^{-4/3}}^{1}
(D_1{g_2})(u,\al)\cR(-\al^{-4/3}X_1,X_1 u) \dd u.
\]
\end{lem}

In the expression for $E_2(\bxi)$ we have used the fact that
$g_2(u,v)=0$ for $u\geq 1$, as follows from \eqref{g2'}.
Much as in the initial calculation of $N(\bxi,\tau_1)$, we shall have
to work rather hard to handle the overall contribution from the error
term $E_2(\bxi)$ in this estimate.  Whereas the error term in Lemma
\ref{lem:t-2-ell} ultimately made a negligible
final contribution to $\#\E(B)$, we
will show that this is no longer the case here.

\begin{pro}\lab{lemma11}
Let $\varepsilon>0$.  Then there exists $\be \in \R$ such that
$$
\sum_{\bxi\in\mcal{A}(B,1)}
E_2(\bxi) =\be B + O_\varepsilon(B^{25/28+\varepsilon}),
$$
where $\mcal{A}(B,1)$ is given by \eqref{defAB}.
\end{pro}

The proof of this result is undertaken in the following section. In
particular the explicit value of $\be$ is given below in (\ref{defbeta}).

\subsection{Proof of Proposition \ref{lemma11}}\lab{error-11}

Throughout this section we shall make frequent use of the inequality
$\log |\bxi| \ll \log B$, that follows from (\ref{h-al}).  We begin
with a rather crude upper bound for
$\cR(t_1,t_2)$, as given by (\ref{bike2}).
Recall the definition of $r(t_1/k_1,t_2/k_1;a\vr^2,q)$ from 
Lemma~\ref{cong}, together with that of $q$ from \eqref{q}. Then 
Lemma~\ref{lem:sum_square} yields
\begin{equation}\lab{majR}
    \begin{split}
\cR(t_1,t_2)
&\ll_\varepsilon
\phi^*(\xi_1\xi_3)\sum_{k_\ell\mid \xi_4\xi_5\xi_6}{|\mu(k_\ell)|\over
k_\ell}\sum_{k_1\mid \xi_2\xi_3\xi_6}|\mu(k_1)| q^{1/2+\varepsilon}\\
&\ll_\varepsilon
2^{\omega(\xi_4\xi_5\xi_6)+\omega(\xi_2\xi_3\xi_6)}(\xi_\ell^3\xi_4^2\xi_5)^{1/2+\varepsilon}\\
&\ll_\varepsilon
4^{\omega(\xi_2\xi_3\xi_6)}(\xi_\ell^3\xi_4^2\xi_5)^{1/2+2\varepsilon},
    \end{split}
\end{equation}
for any $\varepsilon>0$.

Ultimately, our proof of Proposition
\ref{lemma11} will involve an initial summation of
$E_2(\bxi)$ over the variable $\xi_1$.  It will therefore be
convenient to define
$$
Y_1:=\frac{B^{1/2}}{ \xi^{(0,3/2,2,3/2,2,5/2,3)}}=\frac{\xi_1}{\al},
$$
where $\al$ is given by \eqref{notation-1}. In particular this
notation takes the place of that introduced in \eqref{Y12}, and it
it is clear that
$Y_1$ is independent of $\xi_1$.  Moreover, the inequality $\alpha\leq
1$  is plainly equivalent to $\xi_1\leq Y_1$.
We also set
$$
q_1:=\xi^{(0,1,1,2,2,2,2)},\quad
q_2:=\xi^{(0,2,2,1,2,3,4)},
$$
and observe that $X_1={q_1} \alpha^{-2/3}$ in (\ref{notation-2}).
On recalling the definition of $X_2$, we therefore conclude that
$E_2(\bxi)$ can be rewritten as
\begin{equation}\lab{formR'}
\hspace{-0.1cm}
E_2(\bxi)
=
\frac{-B^{1/2}}{ \xi^{(0,1/2,0, 3/2, 1, 1/2,0)}}
\int_{-\alpha^{-4/3}}^{ 1}\!\!\!\! (D_1g_2)(u,\alpha ){\mcal{R}}
(-\alpha^{-2}{q_1},\alpha^{-2/3} {q_1} u)\d u.
\end{equation}

Let $T\geq 1$ be a parameter, to be chosen in due course, and recall
the definition of $\mcal{A}(B,T)$ from \eqref{defAB}.
Our first task is to show that we obtain a satisfactory contribution
by summing $E_2(\bxi)$ over those values of $\bxi$ for which either
$q_2>T$, or $\xi_1 \leq Y_1/T^{1/2}$.

\begin{lem}\lab{except}
Let $\varepsilon>0$.  Then for any $T\geq 1$ we have
$$
\sum_{\colt{\bxi\in \mcal{A}(B,1)}{
\mbox{\scriptsize{$q_2>T$ or $\xi_1\leq Y_1/T^{1/2}$}}
}} E_2(\bxi) \ll_\varepsilon BT^{\ve-1/2}.
$$
\end{lem}

\begin{proof}
On inserting (\ref{majR}) into (\ref{formR'}), and applying the bound
\eqref{intD1g1-b}, we obtain
\begin{equation}\lab{majR'}
E_2(\bxi)  \ll_\varepsilon \frac{B^{1/2}}{
\xi^{(0,1/2,0, 3/2, 1, 1/2,0)}}
4^{\omega(\xi_2\xi_3\xi_6)}(\xi_\ell^3\xi_4^2\xi_5)^{1/2+\varepsilon}.
\end{equation}
It therefore follows that
\begin{align*}
\sum_{\xi_1\leq Y_1 }E_2(\bxi)
&\ll_\varepsilon
\frac{4^{\omega(\xi_2\xi_3\xi_6)}(\xi_\ell^3\xi_4^2\xi_5)^{\varepsilon}
B}{\xi^{(0,2,2, 3/2 , 2, 5/2,3)}}\\
&\ll_\varepsilon
\frac{B}{{q_2}^{1/2-3\ve}\xi^{(0,1+\varepsilon,1+\varepsilon,1+\varepsilon,1+\varepsilon,1+\varepsilon,1+\varepsilon)}}, 
\end{align*}
whence the overall contribution from the case $q_2> T$ is
$O_\varepsilon(BT^{\ve-1/2}),$ on redefining the choice of
$\varepsilon$.  This is satisfactory for the lemma.

In the same fashion one deduces from the upper bound \eqref{majR'} that
\begin{align*}
\sum_{\xi_1\leq Y_1/T^{1/2}}E_2(\bxi)
&\ll_\varepsilon\sum_{\xi_1\leq
Y_1/T^{1/2}}{4^{\omega(\xi_2\xi_3\xi_6)}(\xi_\ell^3\xi_4^2\xi_5)^{1/2+\varepsilon}B^{1/2}
\over \xi^{(0,1/2,0, 3/2, 1, 1/2,0)}}\\
&\ll_\varepsilon
{4^{\omega(\xi_2\xi_3\xi_6)}(\xi_\ell^3\xi_4^2\xi_5)^{\varepsilon}B
\over
T^{1/2}\xi^{(0,2,2, 3/2 , 2, 5/2,3)}},
\end{align*}
and so the overall contribution from the
case $\xi_1\leq Y_1/ T^{1/2}$ is $O(B/T^{1/2}).$ This too is
satisfactory for the lemma.
\end{proof}

It is interesting to remark that on taking $T=1$ in Lemma
\ref{except}, we deduce that the overall contribution from the
$E_2(\bxi)$ term in Lemma \ref{lem:t-1} is $O(B)$. This is already
enough to establish a version of Theorem \ref{main'} with the weaker
error term $O(B)$.

In estimating the overall contribution obtained by summing $E_2(\bxi)$ over
all of the relevant values of $\bxi$, it henceforth suffices to focus our
attention upon those values of $\bxi$ for which
$q_2\leq T$ and $\xi_1 > Y_1/T^{1/2}$.
Let $t_2 \geq t_1$ 
and recall the definition (\ref{bike2})
of $\cR(t_1,t_2)$.  Then it follows that
\begin{align*}
\sum_{\colt{\xi_1\leq Y_1}{\hcf(\xi_1,\xi_2\xi_\ell\xi_4\xi_5)=1}} 
\hspace{-0.6cm}
\cR(t_1,t_2)
&=
\sum_{\colt{k_\ell\mid \xi_4\xi_5\xi_6}{
\hcf(k_\ell, \xi_2\xi_3)=1}}
\hspace{-0.4cm}
{\mu(k_\ell)\over
k_\ell}\sum_{\colt{k_1\mid\xi_2\xi_3\xi_6}{\hcf(k_1,q)=1}}
\hspace{-0.3cm}
\mu(k_1) \sum_{\xi_1\leq Y_1} \varpi (\xi_1)\cR_0(t_1,t_2),
\end{align*}
with
$$
\varpi(\xi_1):=\left\{
\begin{array}{ll}
\phi^*(\xi_1\xi_3), & \mbox{if $\hcf(\xi_1,q\xi_2)=1$},\\
0, & \mbox{otherwise},
\end{array}
\right.
$$
and
$$
\cR_0 (t_1,t_2):= \sum_{\colt{1\leq \vr\leq q}{\hcf(\vr,q)=1}}
r(t_1/k_1,t_2/k_1;a\vr^2,q).
$$
%t%MODIFICATION:
The function $\mcal{R}_0(t_1,t_2)$ is differentiable with respect to
$t_j$ everywhere outside the discrete set $E$ of 
integers congruent to $k_1a \vr^2$ modulo ${q}$, for
some integer $\vr$ which is coprime to $q$.
It will be useful to record the equalities
\begin{equation} \lab{partialR'b1b2}
\frac{\partial  \cR_0}{\partial t_j}(t_1,t_2) =(-1)^{j-1}
\frac{\phi^*(q)}{k_1}, \quad (j=1,2),
\end{equation}
that is valid outside $E$.  Moreover, we shall make use of the bound
\begin{equation} \lab{majR'b1b2}\cR_0 (t_1,t_2)
\ll_\varepsilon q^{1/2+\varepsilon},
\end{equation}
that follows immediately from an application of Lemma \ref{lem:sum_square}.
%t%TO HERE

Returning to the above calculation, we may clearly combine it
with \eqref{formR'} to deduce that
\begin{align}
\begin{split}
\sum_{\colt{Y_1/T^{1/2}< \xi_1\leq Y_1}
{\hcf(\xi_1,\xi_2\xi_\ell\xi_4\xi_5)=1} }
\hspace{-0.3cm}
E_2(\bxi)
&={-B^{1/2}\over
\xi^{(0,1/2,0, 3/2, 1, 1/2,0)}}
\hspace{-0.2cm}
\sum_{\colt{k_\ell\mid
\xi_4\xi_5\xi_6}{\hcf(k_\ell, \xi_2\xi_3)=1}}
\hspace{-0.2cm}
{\mu(k_\ell)\over
k_\ell}
\hspace{-0.2cm}
\sum_{\colt{k_1\mid
\xi_2\xi_3\xi_6}{\hcf(k_1,q)=1}}
\hspace{-0.2cm}
\mu(k_1)\\
&\qquad
\sum_{ Y_1/T^{1/2}<\xi_1\leq
Y_1}\!\!  \varpi (\xi_1)G(\xi_1/Y_1),
\end{split}\lab{sumxi1R'}
\end{align}
with
\begin{equation}\lab{defG}
G(\alpha):=
\int_{g_1(\alpha )}^{1} (D_1g_2)(u,\alpha )\cR_0
(-\alpha^{-2}{q_1},\alpha^{-2/3} {q_1} u)\d u.
\end{equation}
In particular it follows from \eqref{intD1g1-b} and
\eqref{majR'b1b2} that
\begin{equation} \lab{majG}
G(\alpha)\ll_\varepsilon q^{1/2+\varepsilon}.
\end{equation}
We shall also need the following complementary result about the size
of $G'$.

\begin{lem}\lab{lemmajG'} Let $\varepsilon>0$ and let $\alpha\in(0,1]$
such that $ \alpha^2\neq 2/(1+\sqrt{5})$. Then
$$
G'(\alpha)\ll_\varepsilon
\frac{q^{1/2+\varepsilon}}{\alpha^2}
+\frac{q_1}{k_1\alpha^{5/3}}.
$$
\end{lem}

\begin{proof}%t%MODIFICATION
Recall the definitions \eqref{g1}, \eqref{g2'} of the functions $g_1$
and $g_2$.  When
$1+1/v^2\neq 1/v^{4}$, it is clear that $g_1(v)$ is differentiable.
Assume therefore that $ \alpha^2\neq 2/(1+\sqrt{5})$, with $\al \in
(0,1]$, and note that the function $(D_2D_1 g_2)(u,v)$ is identically zero.

One can break the integral in \eqref{defG} into a discrete sum
of integrals that avoid all of places where the integrand is not
differentiable.
%t%MORE DETAILS??
Combining \eqref{intD1g1-a} with \eqref{partialR'b1b2}, a little
thought therefore reveals that
\begin{align}%t%I GET DIFFERENT signs
\begin{split}
G'(\alpha)  =&-{g_1'(\alpha ) }
(D_1g_2)(g_1(\alpha ),\alpha )
\cR_0 (-{q_1}\alpha^{-2},  {q_1} g_1(\alpha )
\alpha^{-2/3}) \\
&- \frac{2\phi^*(q){q_1}g_2(g_1(\alpha ),\alpha )}{
    k_1\alpha^{3}}
+\frac{2\phi^*(q) {q_1}}{
3\alpha^{5/3}k_1}\int_{g_1(\alpha )}^1u(D_1g_2)(u,\alpha )\d u.
\end{split}\lab{calculG'}\end{align}
%t%TO HERE
We use \eqref{majR'b1b2} to estimate the first term,  together
with the obvious bounds 
$$
g_1'(\alpha )\ll 1/\alpha^{5/3},
\quad (D_1g_2)(g_1(\alpha ),\alpha )\ll g_1(\alpha )^{1/2}\ll
1/\alpha^{1/3}.
$$
Thus the first term is easily seen to be
$O_\varepsilon\big(q^{1/2+\varepsilon} /\alpha^2\big).$
Next we deduce from \eqref{g=0} that $g_2(g_1(\alpha ),\alpha )\ll
\alpha^3$, since here the left hand side is zero for $\alpha$
sufficiently
small and bounded otherwise. The second term in
\eqref{calculG'} is therefore bounded by $O(q_1/k_1 )$.
Finally we note that it is straightforward to show that  the integral
in the third term is bounded.  Thus the third term in \eqref{calculG'}
is $O\big(q_1/(k_1\alpha^{5/3})\big)$, which completes the proof of Lemma \ref{lemmajG'}.
\end{proof}

We now turn to the problem of estimating
$$
S(Y):=\sum_{\xi_1\leq Y}  \varpi(\xi_1),
$$
for any $Y\geq 1$, in preparation for an application of integration by
parts in \eqref{sumxi1R'}.  Let
$$
\phi'(n):=\prod_{p\mid n}\Big(1+\frac{1}{p}\Big)^{-1},
$$
for any $n \in \N$. Then we shall establish the following result
rather easily.

\begin{lem}\lab{estS}
We have $S(Y)=S_P(Y)+S_R(Y)$ for any $Y\geq 1$, with
$$
S_P(Y)=\frac{6}{ \pi^2}
\phi'(k_\ell\xi_2\xi_3\xi_\ell\xi_4\xi_5)Y,\quad
S_R(Y)=O_\varepsilon(2^{\omega(k_\ell\xi_2
\xi_\ell\xi_4\xi_5)}Y^{\varepsilon}).
$$
\end{lem}

\begin{proof}
To establish the lemma we consider the corresponding Dirichlet series
\begin{align*}
\sum_{\xi_1=1}^\infty\frac{  \varpi (\xi_1)}{
\xi_1^s}&=\phi^*(\xi_3)
\prod_{p\nmid k_\ell \xi_2\xi_3\xi_\ell\xi_4\xi_5}\Big(1+\frac{1-1/p}{
p^s-1}\Big)
\prod_{p\mid \xi_3} \Big(
\frac{1 }{ 1- p^{-s}}\Big)\\
&=\phi^*(\xi_3)\zeta(s)
\prod_{p\nmid k_\ell\xi_2\xi_3\xi_\ell\xi_4\xi_5 }\Big(1-\frac{1 }{
p^{s+1}}\Big)
\prod_{p\mid k_\ell \xi_2\xi_\ell\xi_4\xi_5 } 
   ( 1- p^{-s}),
\end{align*}
where we have used the relation
$\hcf(\xi_3,k_\ell)=1$.  The result is now immediate on applying
Perron's formula.
\end{proof}

We are now in good shape to complete the proof of Proposition \ref{lemma11}.
Recall the definition \eqref{defG} of the function $G$.  Then it is
easily seen that $G$ is piecewise differentiable, apart from at the
points $\sigma$ that are given by ${q_1}\sigma^{-2}/k_1=a\vr^2+kq$.  At
these points $G$ is discontinuous, but is continuous to the left and
to the right of such points. An integration by parts therefore yields
\begin{align}\begin{split}
\sum_{Y_1/T^{1/2}<\xi_1\leq Y_1}
\hspace{-0.3cm}
   \varpi(\xi_1)G(\xi_1/Y_1)
=& \int_{Y_1/T^{1/2}}^{ Y_1}G(t/Y_1)\d S(t)\\
=&S(Y_1)G(1) -S(Y_1/T^{1/2})G(1/T^{1/2})\\
&-\frac{1}{ Y_1}\int_{Y_1/T^{1/2}}^{
Y_1}G'(t/Y_1)S(t)\d t
\\
&-\int_{1/T^{1/2}}^{1}\sum_\sigma
\big(G(\sigma+)-G(\sigma-)\big)S(Y_1t)\delta_\sigma(t),
\end{split}\lab{intpart}
\end{align}
where the sum in the final term is over the
discontinuities $\sigma$ of $G$, and
$\delta_\sigma$ denotes the Dirac measure for
$\sigma$.

We begin by estimating the contribution to (\ref{intpart}) from the
principal term $S_P$ in the estimate for $S$.  This gives
\begin{equation}
\int_{Y_1/T^{1/2}}^{ Y_1}G(t/Y_1)\d
S_P(t)= \frac{6}{
\pi^2}  {Y_1 } \phi'(k_\ell\xi_2\xi_3\xi_\ell\xi_4\xi_5)  \int_{1/T^{1/2}}^1
G(w) \d w.
\lab{TP}\end{equation}
The overall contribution in (\ref{sumxi1R'}) from $S_P$ is therefore
\begin{equation}
\begin{split}
&=B\frac{-6}{ \pi^2}
    \sum_{\colt{k_\ell\mid
\xi_4\xi_5\xi_6}{
\hcf(k_\ell, \xi_2\xi_3)=1}}
\!\!\!\!\!\!\frac{\mu(k_\ell)\phi'(k_\ell\xi_2\xi_3\xi_\ell\xi_4\xi_5)}{
k_\ell \xi^{(0,2,2, 3, 3,
3,3)}}\!\!\!\!\! \sum_{\colt{k_1\mid
\xi_2\xi_3\xi_6}{\hcf(k_1,q)=1}}\!\!\!\!\!\!
\mu(k_1)  \int_{1/T^{1/2}}^1  G(w) \d w\\
&=
Bu(\bxi)+O_\varepsilon(B/T^{1/2-\varepsilon}),
\end{split}\lab{TP2}
\end{equation}
where we have written
$$
u(\bxi):=\frac{-6}{ \pi^2}
    \sum_{\colt{k_\ell\mid
\xi_4\xi_5\xi_6}{
\hcf(k_\ell, \xi_2\xi_3)=1}} 
\!\!\!\!\!\!\frac{\mu(k_\ell)\phi'(k_\ell\xi_2\xi_3\xi_\ell\xi_4\xi_5)}{
k_\ell \xi^{(0,2,2, 3, 3,
3,3)}}\!\!\!\!\! \sum_{\colt{k_1\mid
\xi_2\xi_3\xi_6}{\hcf(k_1,q)=1}}\!\!\!\!\!\!
\mu(k_1)  \int_{0}^1  G(w) \d w.
$$
Here we recall that the function $G$, as given by \eqref{defG},
depends intimately upon the values of $k_1$ and $k_\ell$.
The quantity $u(\bxi)$ is the general term of an absolutely convergent
series which is completely independent of $B$.  In fact, by using
\eqref{majG} it is easy to establish the upper bound
$$
u(\bxi)\ll_\varepsilon
2^{\omega(\xi_2\xi_3\xi_6)+\omega(\xi_4\xi_5\xi_6)}\frac{(
\xi_\ell^3\xi_4^2\xi_5)^{1/2+\varepsilon}}{  \xi^{(0,2,2, 3, 3,
3,3)}}\ll_\varepsilon {1\over  {q_2}^{1/2-3\varepsilon}\xi^{ (0,1,1, 1, 1,
1,1)}},
$$
whence
$$
\sum_{\colt{\xi_2,\ldots,\xi_6 \in \N}{q_2>T}}|u(\bxi)|\ll_\varepsilon {1\over
T^{1/2-4\varepsilon}}.
$$
On summing (\ref{TP2}) over all relevant values of $\xi_2,\ldots,\xi_6$
such that $q_2\leq T$, we therefore obtain the overall contribution
\begin{equation}\lab{C1}=\beta B+O\Big({B\over
T^{1/2-4\varepsilon}}\Big)\end{equation}
to (\ref{sumxi1R'}) from $S_P$, where
\begin{equation}
\lab{defbeta}
\beta:=\sum_{\colt{\xi_2,\ldots,\xi_6\in
\N}{\hcf (\xi_\ell,
\xi_2\xi_3)=1}}|\mu(\xi_2\xi_3\xi_4\xi_5)|u(\bxi).
\end{equation}

We proceed to estimate the overall contribution to the second line in
(\ref{intpart}) from the residual term $S_R$ in Lemma \ref{estS}.  We shall need the observation that
$$
\int_{1/T^{1/2}}^{1}|G'(\alpha)|
\d \alpha\ll_\varepsilon
q^{1/2+\varepsilon}T^{1/2} +{q_1} T^{1/3}.
$$
that follows easily from Lemma \ref{lemmajG'}.
Using this bound therefore yields the overall contribution
$$
\ll_\varepsilon
Y_1^{\varepsilon}2^{\omega(\xi_2\xi_3\xi_6)+\omega(\xi_4\xi_5\xi_6)}
\big( q^{1/2+\varepsilon}T^{1/2}+ {q_1}
T^{1/3}\big)
\ll_\varepsilon B^\varepsilon\big(
q^{1/2}T^{1/2} +{q_1} T^{1/3}\big),
$$
to the second line.  Given that $q_2\leq T $, the contribution
in (\ref{sumxi1R'}) is therefore
\begin{align*}
&\ll_\varepsilon
{B^{1/2+\varepsilon}\over
\xi^{(0,1/2,0, 3/2, 1,
1/2,0)}}\Big((\xi_\ell^3\xi_4^2\xi_5)^{1/2}T^{1/2}\frac T{q_2} +
T^{11/6}\frac{q_1}{{q_2}^{3/2}}\Big)
\ll_\varepsilon   {B^{1/2+\varepsilon} T^{11/6}
\over
\xi^{(0,5/2,2, 1, 2, 3,4)}},
\end{align*}
whence we obtain a total contribution of
\begin{equation}\lab{C2}
O( B^{1/2+ 2\varepsilon} T^{11/6}).
\end{equation}

It remains to the examine the influence of the discontinuities of $G$
upon the final result.  In other words, we must now estimate the
overall contribution from the final line (\ref{intpart})
when $S$ is replaced by $S_R$.
Now we have
$$
\sum_\sigma
\big(G(\sigma+)-G(\sigma-)\big)
\delta_\sigma(\alpha) =  \sum_{\colt{1\leq \vr\leq q}{\hcf(\vr,q)=1}}
\delta_{-a\vr^2+q\Z}( {q_1}/\alpha^2k_1)g_2(g_1(\alpha ),\alpha ).
$$
Let $\alpha_0=1/2^{1/2}$.  Then in view of \eqref{g=0} we have
$g_2(g_1(\alpha ),\alpha )=0$  for any $\alpha\leq \alpha_0$, and so
for each $k_1$ it suffices to consider  the values of
${q_1}/\alpha^2k_1$ in $ -a\vr^2+q\Z$ such that $\alpha>\alpha_0$. Hence
\begin{align*}
\sum_\sigma
\big(G(\sigma+)-G(\sigma-)\big)
\delta_\sigma(\alpha)&=
    \sum_{\colt{m\leq {q_1}/\alpha_0^2k_1}{\hcf(m,q)=1}}
\hspace{-0.1cm}
\sum_{\colt{1\leq \vr\leq q}{-a\vr^2\equiv m\mod q}}
\hspace{-0.6cm}
\delta_{m}( {q_1}/\alpha^2k_1)g_2(g_1(\alpha ),\alpha ) .
\end{align*}
Now it is clear that there are at most $2^{\omega(q)}$
solutions $\vr$ to the congruence $-a\vr^2\eqm{m}{q}$, whence
$$
\int_{1/T^{1/2}}^{1}\sum_\sigma
\big(G(\sigma+)-G(\sigma-)\big)S_R(Y_1t)\delta_\sigma(t)
\ll_\varepsilon  {q_1} B^{\varepsilon}.
$$
This therefore yields a contribution of
$$
\ll_\varepsilon
{B^{1/2+\varepsilon}{q_1T^{3/2}}\over
\xi^{(0,1/2,0, 3/2, 1, 1/2,0)}{q_2}^{3/2}} \ll_\varepsilon
{B^{1/2+\varepsilon}T^{3/2}\over
\xi^{(0,5/2,2, 1, 2, 3,4)}}
$$
in (\ref{sumxi1R'}), since $q_2\leq T$.
Thus the  total contribution is
\begin{equation}\lab{C3}
O_\varepsilon(B^{1/2+2\varepsilon}T^{3/2}).
\end{equation}

On combining the different contributions that we have estimated in
(\ref{C1}), (\ref{C2}) and (\ref{C3}), and then combining this with
Lemma \ref{except}, we have therefore established the estimate
$$
\sum_{\bxi\in\mcal{A}(B,1)} E_2(\bxi)
=\beta B+
O_\varepsilon\Big(\frac{B}{T^{1/2-4\varepsilon}}+B^{1/2+2\varepsilon} T^{11/6}
\Big),
$$
for any $T \geq 1$.  Making the choice $T=B^{3/14}$
therefore allows us to complete the proof of Proposition
\ref{lemma11}, on redefining the choice of $\varepsilon$.

\subsection{Summation over $\bxi$}\lab{1-6}

The aim of this short section is to sum the main term
in Lemma~\ref{lem:t-1}'s estimate for $N'(\bxi)$ over all
$\bxi \in \FF$ such that the height condition \eqref{h-al} holds.
Note that the definition \eqref{vt} of $\vt$ ensures that this main
term is zero unless $\bxi\in\FF$.

Define the arithmetic function
\begin{equation}\lab{D}
\Delta(n) = B^{-5/6}\sum_{\colt{\bxi\in \N^7}{\base 2 3 4 3 4 5 6 = n}}
\frac{\vartheta(\bxi) X_1X_2}{\xi_\ell^3\xi_4^2\xi_5},
\end{equation}
for any $n \in\N$,
where $X_1,X_2$ are given by (\ref{notation-2}).
Let $\ve>0$ and let $B \geq 1$ be such that
$B\in\R\smallsetminus\overline\Q$.
Then on bringing together the estimates in Lemmas
\ref{lem:t-2-ell} and \ref{lem:t-1}, in addition
to the
handling of the error terms in Propositions \ref{lemma10} and
\ref{lemma11}, we therefore deduce that
$$
\#\E(B)=
B^{5/6}\sum_{n \leq B} \Delta(n) g_3((n/B)^{1/6}) + \be B +
O_\varepsilon(B^{43/48+\varepsilon})
$$
in Lemma \ref{lem:count}.
But we may now combine this with Lemmas \ref{Reduc1} and
\ref{lem:count} to deduce the following basic result.

\begin{lem}\label{lem:NUB}
Let $\varepsilon>0$ and let $B \geq 1$.  Then we have
\[
\NU(B)=2B^{5/6}\sum_{n \leq B} \Delta(n)
g_3((n/B)^{1/6}) +
\Big(\frac{12}{\pi^2}+2\be\Big)B +
O_\varepsilon(B^{43/48+\varepsilon}),
\]
where
$g_3$ is given by \eqref{g3},
$\be$ is given by \eqref{defbeta} and
$\D$ is given by \eqref{D}.
\end{lem}

\begin{proof}
The result is immediate if $B\in\R \smallsetminus\overline\Q$. Let
$B_1,B_2, B_3, \ldots$ be an arbitrary sequence of transcendental numbers
such that
$\lfloor B_k \rfloor=\lfloor B\rfloor$ for each $k \in \N$ and $\lim_{k\to \infty}B_k=B$.
Then for each $k \in \N$ we have
\begin{align*}
\NU(B)&=\NU(B_k)\\&=2B^{5/6}\sum_{n \leq B} \Delta(n) g_3((n/B_k)^{1/6}) +
\Big(\frac{12}{\pi^2}+2\be\Big)B + O_\varepsilon(B^{43/48+\varepsilon}).
\end{align*}
Taking the limit as $k\to \infty$ therefore yields the statement of Lemma \ref{lem:NUB}.
\end{proof}

\section{Proof of Theorem \ref{main'}: d\'enouement}\lab{proof4'}

In this section we complete the proof of Theorem \ref{main'}.
Our main task is to estimate the quantity
\[
M(x) := \sum_{n \le x}
\Delta(n),
\]
for any $x \geq 1$.   This will be achieved in Lemma \ref{chalk} below.
First it will be necessary to examine the
analytic properties of the corresponding Dirichlet series
$F(s):=\sum_{n=1}^\infty \D(n)n^{-s}$.
On recalling the definition \eqref{D} of $\D$,
and those of $X_1,X_2$, we see that
\[
\Delta(n) = \sum_{\colt{\bxi\in \N^7}{\base 2 3 4
3 4 5 6 = n}} \frac{\vartheta(\bxi)
n^{1/6}}{\base 1 1 1 1 1 1 1}
\]
in \eqref{D}. Hence
\begin{equation}\lab{egg}
F(s-5/6)=\sum_{\bxi\in \N^7}
\frac{\vartheta(\bxi)\base 2 3 4 3 4 5
 6}{\xi_1^{2s+1}\xi_2^{3s+1}\xi_3^{4s+1}\xi_\ell^{3s+1}\xi_4^{4s+1}\xi_5^{5s+1}\xi_6^{6s+1}}, 
\end{equation}
where $\vartheta(\bxi)$ is given by \eqref{vt}.
Now recall the definitions \eqref{E1}, \eqref{E2} of $E_1(s)$ and
$E_2(s)$,
and define the half-plane
\begin{equation}\lab{H-plane}
\HH_{\theta} := \{s \in \C: ~{\Re e}(s) > \theta \}
\end{equation}
for any $\theta>0.$  We proceed by establishing the following result.

\begin{lem}\label{F}
Let $\varepsilon>0$. Then there exists  a function $G_{1,1}(s)$ that
is holomorphic and bounded function on $\HH_{5/6+\varepsilon}$,
such that
$$
F(s-5/6)=E_1(s)E_2(s)G_{1,1}(s).
$$
\end{lem}

\begin{proof}
We have already seen in \eqref{egg} an explicit formula for
$F(s)$.  On writing $F(s+1/6) = \prod_p F_p(s+1/6)$ as a product of
local factors, a straightforward calculation reveals that
\begin{equation*}
\begin{split}
F_p(s+1/6) = &1 + \frac{(1-1/p)^2}{(p^{6s+1}-1)}\Big(
     \frac{p^{2s+1}}{p^{2s+1}-1} +
     \frac{p^{2s+1}p^{6s+1}}{p^{4s+1}(p^{2s+1}-1)}\\
     &+
     \frac{p^{6s+1}}{(1-1/p)p^{3s+1}} +  \frac{1}{p^{3
s+1}-1}+ \frac{p^{3 s+1}p^{6s+1}}
     {p^{4s+1}(p^{3 s+1}-1)} \\
     &+ \frac{p^{3 s+1}
       p^{6s+1}}{p^{5s+1}(p^{3 s+1}-1)}\Big) +
     \frac{1-1/p}{p^{2s+1}-1} + \frac{1-1/p}{p^{3 s+1}-1}.
\end{split}
\end{equation*}
We proceed to calculate
$F_p(s+1/6)(1-1/p^{6s+1})(1-1/p^{2s+1})(1-1/p^{3s+1})$
on
$\HH_{-1/6+\varepsilon}$. Thus we have
\begin{equation*}
\begin{split}
F_p(s+1/6)&(1-1/p^{6s+1})(1-1/p^{2s+1})(1-1/p^{3s+1}) \\
&=  1 + \frac 1{p^{3s+1}} +  \frac 2{p^{4s+1}}+  \frac 1{p^{5s+1}}- \frac
2{p^{7s+2}}  + O_\varepsilon\Big(\frac 1 {p^{1+\varepsilon}}\Big),
     \end{split}
\end{equation*}
whence
\begin{equation*}
\begin{split}
       \frac{F_p(s+1/6)}{E_{1,p}(s+1)}  = &1 - \frac 1 {p^{10s+2}}- \frac 2
{p^{9s+2}} - \frac 4{p^{8s+2}} - \frac 4 {p^{7s+2}} \\& + \frac 2 {p^{14s+3}}+
\frac 5 {p^{13s+3}} + O_\varepsilon\Big(\frac 1 {p^{1+\varepsilon}}\Big)\\
          =&
E_{2,p}(s+1)\Big(1+O_\varepsilon\Big(\frac 1
{p^{1+\varepsilon}}\Big)\Big).
\end{split}
\end{equation*}
This therefore establishes the lemma.
\end{proof}

Let $\varepsilon>0$ and recall the properties of $E_1(s), E_2(s)$ that were
outlined in \S \ref{intro}.  In particular we clearly have
$$
E_1(s+1)=\prod_{i\in\{1,2,3,\ell,4,5,6\}}
\zeta(\la_i s+1),
$$
where $\bla=(\la_1,\ldots,\la_6)$ is given by \eqref{bla},
and so
$$
E_1(s)=\frac{1}{\prod_i \la_i}(s-1)^{-7}+O\big((s-1 )^{-6}\big),
$$
as $s\rightarrow 1$.
It follows from Lemma \ref{F} that
$F(s-5/6)$ is meromorphic on $\HH_{9/10+\varepsilon}$, with a pole of
order $7$ at $s=1$.  Moreover, the function
$$
G(s):=\frac{F(s-5/6)}{E_1(s)}
$$
is holomorphic and bounded function on $\HH_{9/10+\varepsilon}$.
Now a simple calculation reveals that for any $x\geq 1$ we have
\begin{equation}\lab{res}
\rom{Res}_{s=1}\Big(\frac{F(s-5/6)x^{s-5/6}}{s-5/6}\Big)=
\frac{6 G(1)x^{1/6}Q_0(\log x)}{6! \cdot \prod_i \lambda_i},
\end{equation}
for some monic polynomial $Q_0$ of degree $6$,
and where
\begin{equation}\lab{G(1)}
G(1)=
\prod_p\Big(1-\frac{1}{p}\Big)^7\Big(1+\frac 7 p + \frac 1 {p^2}\Big).
\end{equation}
We are now ready to establish the following
estimate for $M(x)$.

\begin{lem}\lab{chalk}
Let $\varepsilon>0$.  Then there exists a monic polynomial
$Q_1$ of degree $6$ such that
$$
M(x) =
\frac{G(1)x^{1/6}Q_1(\log x)}{1036800} +
O_\varepsilon(x^{1/6-1/11+\varepsilon}),
$$
for any $x \geq 1$.
\end{lem}

\begin{proof}
Our starting point is the pair of inequalities
\begin{equation}\lab{viva}
{1\over y}\int_{x-y}^xM(t)\d t\leq M(x)\leq {1\over y}\int_{x
}^{x+y}M(t)\d t,
\end{equation}
that are valid for any $1\leq y<x$.
Note that $M(x)$ is an increasing function. This approach to
estimating $M(x)$  is a variant of Perron's formula which has the advantage
of giving rise to absolutely convergent complex integrals.
We shall only estimate the right-hand side of \eqref{viva}, since the integral
on the left ultimately yields precisely the same
estimate, as is easily checked.
An application of Perron's formula yields
\begin{align*}
{1\over y}\int_{x }^{x+y}M(t)\d t
&= {1\over 2\pi i}\int_{1/6+1/\log x-i\infty}^{1/6+1/\log x+i\infty}
F(s){(x+y)^{s+1}-x^{s+1}\over y s(s+1)} \d s\\
&= {1\over 2\pi i}\int_{\kappa-i\infty}^{\kappa+i\infty}
F(s-5/6){(x+y)^{s+1/6}-x^{s+1/6}\over y (s-5/6)(s+1/6)} \d s,
\end{align*}
with $\kappa=1+1/\log x$.
Let $\delta=1/10+\varepsilon$.
Then we shall apply Cauchy's residue theorem to the rectangular
contour $\mcal{C}$ joining the points
$$
\kappa-iT,\quad \kappa+iT, \quad 1-\delta+iT, \quad 1-\delta-iT,
$$
for some value of $T \in [1,x]$ to be selected in due course.

We shall employ the
well-known convexity bounds
$$
\zeta(\sigma+i \tau)\ll_\varepsilon \left\{
\begin{array}{ll}
|\tau|^{(3-4\sigma)/6+\varepsilon},& \mbox{if $\sigma\in [0,1/2]$,}\\
|\tau|^{(1-\sigma)/3+\varepsilon},& \mbox{if $\sigma\in
    [1/2,1]$,}
\end{array}
\right.
$$
valid for any $|\tau|\geq 1$
(see Tenenbaum \cite[\S II.3.4]{ten}, for example).
It follows from these that
$$
F(\sigma-5/6+i\tau)\ll
E_1(\sigma+i\tau)\ll_\varepsilon
(1+|\tau|)^{\mu_F(\sigma)+\varepsilon}
$$
on $\mcal{C}$, with
$$
\mu_F(\sigma)\leq \left\{
\begin{array}{ll}
    \max\{9(1-\sigma),0\}& \mbox{if $\sigma\geq 11/12$},\\
65/6-11\sigma& \mbox{if $9/10\leq \sigma\leq 11/12$}.
\end{array}
\right.
$$
Moreover, we have
\begin{equation}\lab{maj}
(x+y)^{s+1/6}-x^{s+1/6}\ll
x^{\sigma+1/6}(y(1+|\tau|)/x)^\alpha,
\end{equation}
for any $\alpha\in [0,1]$.  On taking $\alpha=1$ here we
may therefore estimate the contribution from the horizontal strips as being
\begin{align*}
\int_{1-\delta\pm iT}^{\kappa\pm iT}
\Big|F(s-5/6) {(x+y)^{s+1/6}-x^{s+1/6}\over y (s-5/6)(s+1/6)}
\Big| \d  s
&\ll_\varepsilon \int_{1-\delta }^{\kappa }
x^{\sigma-5/6}T^{\mu_F(\sigma)-1+\varepsilon}
\d \sigma\\
& \ll_\varepsilon
{x^{1/6+\varepsilon} \over T}+ x^{1/6-\delta +\varepsilon},
\end{align*}
since $T \geq 1$.
In a similar fashion we see that
\begin{align*}
\int_{\kappa\pm iT}^{\kappa\pm i\infty}
\Big|F(s-5/6) {(x+y)^{s+1/6}-x^{s+1/6}\over y (s-5/6)(s+1/6)}
\Big| \d  s
& \ll_\ve  {x^{1/6+\varepsilon} \over T}.
\end{align*}
On applying \eqref{maj}
with $\alpha=3/2-6\delta$,  we see that
the contribution from the line $\Re e (s)=1-\delta$ is
\begin{align*}
&\ll
\int_{1-\delta-iT}^{1-\delta+iT}
\Big|F(s-5/6) {(x+y)^{s+1/6}-x^{s+1/6}\over y (s-5/6)(s+1/6)}
\Big| \d  s \\
&\ll
x^{5\delta-1/3}y^{1/2-6\delta}\int_{-T}^{
T}\Big|{E_1(1-\delta+i\tau)
\over (1+|\tau|)^{1/2+6\delta}}\Big|\d \tau
\\&\ll
x^{5\delta-1/3}y^{1/2-6\delta}\log T  \sup_{1\leq
U\leq T} {1\over U^{1/2+6\delta}}\int_{0}^{  U}\left|{ E_1(1-\delta+i\tau) }
\right|\d \tau.
\end{align*}
Now let
$$
J_k(\sigma,U):=\Big({1\over U}\int_{0}^{   U}\left|
\zeta(1-\sigma+i\tau)
\right|^k\d \tau\Big)^{1/k}.
$$
Then an application of H\"older's inequality, with weights $(10, 5,
10/3,4, 20/3)$, yields
\begin{align*}
{1\over U} \int_{0}^{  U}\left| E_1(1-\delta+i\tau)  \right|\d
\tau
\ll& J_{ 10}( 2\delta,2U)J_{ 10}(3\delta,3U)^2 J_{20/3}(4\delta,4U)^2
\\&\times J_{4}(5\delta,5U) J_{20/3}(6\delta,6U).
\end{align*}
It therefore follows from classical upper bounds for the fractional
moment of the Riemann zeta function (see Ivi\'c \cite[Theorem 2.4]{MR792089},
for example) that
$$
\begin{array}{c}
J_{ 10}( 2\delta,2U)\ll_\varepsilon
U^\varepsilon, \quad J_{10}(3\delta,3U)\ll_\varepsilon U^\varepsilon,
\quad J_{20/3}(4\delta,4U) \ll_\varepsilon
U^\varepsilon,\\
J_{4}(5\delta,5U)\ll_\varepsilon
U^\varepsilon, \quad J_{20/3}(6\delta,6U) \ll_\varepsilon U^{6\delta-1/2},
\end{array}
$$
that are valid when $\varepsilon$ is chosen to be sufficiently small.
To obtain the final estimate here we have used the functional equation
of the Riemann zeta function in the form
$
|\zeta(1-6\delta+i\tau)|\asymp |\tau|^{6\delta-1/2} |\zeta(
6\delta+i\tau)|,
$
for any $|\tau|\geq 1$.
Putting this together we therefore obtain the estimate
\begin{align*}
\int_{1-\delta-iT}^{1-\delta+iT}
\Big|F(s-5/6) {(x+y)^{s+1/6}-x^{s+1/6}\over y (s-5/6)(s+1/6)}
\Big| \d  s
&\ll_\varepsilon x^{5\delta-1/3}y^{1/2-6\delta}T^\varepsilon.
\end{align*}

We may now collect together all of our various estimates, together
with the choices $T=x$ and $y=x^{10/11}$, in order to deduce
from \eqref{bla}, \eqref{res} and \eqref{viva}
that there exists a monic polynomial
$Q_1$ of degree $6$ such that
\begin{align*}
M(x)&=
\rom{Res}_{s=1}\Big(F(s-5/6){(x+y)^{s+1/6}-x^{s+1/6}\over
y (s-5/6)(s+1/6)}\Big)
+O_\varepsilon(x^{1/6-1/11+\varepsilon})\\
&=
\frac{6 G(1)x^{1/6}Q_1(\log x)}{6! \cdot \prod_i \lambda_i} +
O_\varepsilon(x^{1/6-1/11+\varepsilon}),
\end{align*}
for any $x\geq 1$.  This completes the proof of Lemma \ref{chalk},
since $\prod_i \la_i=1036800$ by \eqref{bla}.
\end{proof}

We are now ready to complete the proof of Theorem \ref{main'} via an
application of partial summation in Lemma \ref{lem:NUB}.  Thus it follows
from Lemma \ref{chalk} that
there exists a monic polynomial $Q_2$ of degree $6$ such that
\begin{equation*}
    \begin{split}
\sum_{n \leq B} \Delta(n) g_3((n/B)^{1/6})=&
\frac{G(1)}{1036800}
\int_0^B g_3((n/B)^{1/6})
\frac \dd{\dd u}(u^{1/6}Q_1(\log u)) \dd u\\
&+ O_\varepsilon(B^{1/6-1/11+\varepsilon})\\
=&
\frac{G(1)B^{1/6}Q_2(\log B)}{1036800}
\int_0^1g_3(v)  \dd v + O_\varepsilon(B^{5/66+\varepsilon}),
    \end{split}
\end{equation*}
for any $B\geq 1$.  On noting that
   \[
\omega_\infty =
12\int_0^1 g_3(v)\d v
\]
in \eqref{om-inf}, and recalling the equality \eqref{G(1)} for
$G(1)$, we may therefore insert this estimate into Lemma~\ref{lem:NUB} in order to
conclude the proof of Theorem~\ref{main'}.

\section{Proof of Theorem \ref{thm:main}}\label{proof4}

The proof of Theorem \ref{thm:main} is almost identical to the
corresponding proof of \cite[Theorem 1]{math.NT/0412086}, and
so we shall be very brief.
Recall the definition \eqref{H-plane} of the half-plane
$\HH_{\theta}$ for any $\theta>0.$
Following \cite[\S 6]{math.NT/0412086} therefore, we easily deduce
from Lemma \ref{lem:NUB} that
$$
Z_{U,H}(s):=\sum_{x \in U(\Q)}\frac{1}{H(x)^s}
=Z_1(s)+\frac{{12/\pi^2}+2\be}{s-1}+G_2(s),
$$
for $s \in \HH_1$, where
\begin{equation}\lab{G2}
G_2(s)=s\int_{1}^\infty t^{-s-1}R(t)\d t
\end{equation}
for some function $R(t)$ such that $R(t)\ll_\varepsilon
t^{43/48+\varepsilon}$,
and
\begin{align*}
Z_1(s)&=
2s\int_1^\infty t^{-s-1/6}\sum_{n\leq t}\Delta(n)g_3
\Big(\Big(\frac{n}{t}\Big)^{1/6}\Big)\d t.
\end{align*}
Now it is clear that $G_2(s)$ is holomorphic on
$\HH_{43/48+\varepsilon},$ and satisfies
$G_2(s)\ll 1+|\Im m (s)|$ on this domain.
Moreover one readily deduces from
\eqref{E1}, \eqref{E2} and Lemma \ref{F} that
\[
Z_1(s) = 2sF(s-5/6)\int_1^\infty t^{-s-1/6}g_3(1/t^{1/6})\d  t =
E_1(s)E_2(s)G_{1,1}(s)G_{1,2}(s)
\]
on $\HH_{1}$.
Here the function $G_{1,1}(s)$
is holomorphic and bounded on $\HH_{5/6+\varepsilon}$, and
\begin{equation}
\label{G12}
G_{1,2}(s) =12s\int_0^1 v^{6s-6}g_3(v)\d v.
\end{equation}
A simple calculation reveals that
$G_{1,2}(1) =\omega_\infty$, in the notation of \eqref{om-inf}.
Moreover, an application of partial integration yields
$$
G_{1,2}(s) =\frac{12s}{6s-5}\Big(g_3(1)-\int_0^1
v^{6s-5}{g_3}'(v)\d v\Big),
$$
whence $G_{1,2}(s)$ is holomorphic and bounded
on $\HH_{5/6+\varepsilon}$.
In view of the fact that  $G_{1,1}(1)G_{1,2}(1) \neq 0$, we may therefore set
\begin{equation}\lab{G1}
G_1(s):=G_{1,1}(s)G_{1,2}(s)
\end{equation}
on $\HH_{5/6+\varepsilon}$, in order to complete
the proof of Theorem~\ref{thm:main}.

\end{document}